\documentclass[reqno,11pt]{amsart}
\usepackage{enumerate}
\usepackage{amssymb}
\usepackage[headings]{fullpage}

\theoremstyle{plain}
\newtheorem{thm}{Theorem}[section]
\newtheorem{prop}[thm]{Proposition} 
\newtheorem{lem}[thm]{Lemma} 
\newtheorem{cor}[thm]{Corollary} 
\theoremstyle{definition}
\newtheorem{defn}[thm]{Definition} 
\newtheorem*{ack}{\textbf{Acknowledgement}}
\theoremstyle{remark}
\newtheorem{rem}[thm]{\textsl{Remark}} 
\newtheorem{exmp}[thm]{\textsl{Example}} 
\numberwithin{equation}{section}

\DeclareMathAlphabet{\mathsfsl}{OT1}{cmss}{m}{sl}
\newcommand\ad{\operatorname{ad}}
\newcommand\ads{\operatorname{ad}^\star}
\newcommand\appli[5]{
   \begin{matrix}#1\colon\vphantom{#2#3}\\*\vphantom{#4#5}\end{matrix}
   \begin{matrix}#2\vphantom{#1#3}\\*#4\vphantom{#5}\end{matrix}\;
   \begin{matrix}\longrightarrow\vphantom{#1#2#3}\\*
     \longmapsto\vphantom{#4#5}\end{matrix}\;
   \begin{matrix}#3\vphantom{#1#2}\\*#5\vphantom{#4}\end{matrix}}
\newcommand\can{\mathit{can}}
\renewcommand\cosh{\operatorname{ch}}
\renewcommand\coth{\operatorname{cotanh}}
\newcommand\cprime{$'$}
\renewcommand\d{\operatorname{d}}

\newcommand\dlie{\lie{d}}

\newcommand\ensemble[2]{\left.\left\{#1\vphantom{#2}
  \,\right\vert\,#2\right\}}
\newcommand\g{\lie{g}}
\newcommand\h{\lie{h}}
\newcommand\happy[1]{{\overset{\smile}{#1}}_p}
\newcommand\ie{i.e.,\ }

\newcommand\klie{\lie{k}}
\renewcommand\l{\lie{l}}
\newcommand\lcan{l^\can}
\newcommand\lcanG[1]{l^{\can(#1)}}
\newcommand\lie[1]{\mathfrak{#1}}
\newcommand\lEV{\ell^\mathit{EV}}
\newcommand\m{\lie{m}}
\newcommand\op{\mathsf{op}}
\newcommand\p{\mathrm{p}}
\newcommand\rAM{r^\mathit{AM}}

\newcommand\rEV{r^\mathit{EV}}
\newcommand\rmx{\mathbf{r}}
\newcommand\saddy[1]{{\overset{\frown}{#1}}_p}
\newcommand\saddyk[1]{{\overset{\frown}{#1}}_\bullet}
\renewcommand\sinh{\operatorname{sh}}
\renewcommand\t{\mathsfsl{t\mskip.2\thinmuskip}}

\newcommand\upsi{{\underline{\psi}}}

\newcommand\ulie{\lie{u}}
\newcommand\A{\mathcal{A}}
\newcommand\Ad{\operatorname{Ad}}

\newcommand\Alg{\operatorname{\mathsfsl{A}}}
\newcommand\Alt{\mathsf{Alt}}
\newcommand\Cycl{\mathop{\circlearrowleft}\limits}
\newcommand\CDYB{\mathsfsl{CDYB}}
\newcommand\DD{\mathbb{D}}
\newcommand\Dynl{\mathsfsl{Dynl}}

\newcommand\G{\mathcal{G}}
\newcommand\GG{\mathbb{G}}
\newcommand\Id{\mathrm{Id}}

\newcommand\Ker{\operatorname{\mathsf{Ker}}}

\renewcommand\L{\mathcal{L}}
\newcommand\Lie{\operatorname{\mathsfsl{Lie}}}
\newcommand\LL{\mathbb{L}}

\newcommand\Map{\mathsfsl{Map}}

\newcommand\Retire[1]{\relax}
\newcommand\Sym{\mathrm{S}}

\newcommand\T{\operatorname{\mathsf{T}}}

\newcommand\<{\mathopen\langle}
\renewcommand\>{\mathclose\rangle}

\def\middletwolines{\\*}
\def\twolinesleft#1#2\twolinesmid#3\twolinesright#4{%
  \left#1\vphantom{#3}#2\right.\middletwolines\left.\vphantom{#2}#3\right#4}

\def\ch{\operatorname{ch}}
\def\sh{\operatorname{sh}}
\def\chpmname{\operatorname{c}}
\def\shpmname{\operatorname{s}}
\def\chx#1{{\chpmname^{{\scriptscriptstyle#1}}}}
\def\shx#1{{\shpmname^{{\scriptscriptstyle#1}}}}
\def\chp{\chx+}
\def\chm{\chx-}
\def\chpm{\chx\pm}
\def\shp{\shx+}
\def\shm{\shx-}
\def\shpm{\shx\pm}

\newcommand\e{\operatorname{e}}

\newcommand\Nset{\mathbb{N}}

\begin{document}
\setlength{\baselineskip}{1.05\baselineskip} 
\title{Bidynamical Poisson groupoids}
\author[R.~Pujol]{Romaric Pujol}\email{pujol@math.toronto.edu}
\address{University of Toronto, Department of Mathematics, 40 St George
Street, Toronto, \mbox{Ontario M4S\ 2E4}, Canada.}
\date{\today}
\begin{abstract}
We give relations between dynamical Poisson groupoids, classical
dynamical Yang--Baxter equations and Lie quasi-bialgebras. We show that
there is a correspondance between the class of bidynamical Lie
quasi-bialgebras and the class of bidynamical Poisson groupoids. We give
an explicit, analytical and canonical equivariant solution of the
classical dynamical Yang--Baxter equation (classical dynamical
$\ell$-matrices) when there exists a reductive decomposition
$\g=\l\oplus\m$, and show that any other equivariant solution is
formally gauge equivalent to the canonical one. We also describe the
dual of the associated Poisson groupoid, and obtain the characterization
that a dynamical Poisson groupoid has a dynamical dual if and only if
there exists a reductive decomposition $\g=\l\oplus\m$.
\end{abstract}
\maketitle
\thispagestyle{empty}

\section*{Introduction}
The Classical Dynamical Yang--Baxter equation (CDYBE) is an important
differential equation in mathematical physics. It was first introduced
by Felder~\cite{F95} in the context of conformal field theory, appearing
as a dynamical analogue of the Classical Yang--Baxter equation (CYBE),
which plays a central role in the theory of integrable systems; the
geometric meaning of (CYBE) was given by Drinfel\cprime{}d, and gives
rise to the theory of Poisson--Lie groups. The geometric meaning of
(equivariant solutions of) the (CDYBE) was given by Etingof and
Varchenko~\cite{EV98}, and is a groupoidal analogue of that of (CYBE):
dynamical Poisson groupoids. In the present paper, we explicitely
describe dynamical Poisson groupoids with base space containing $0$
which have a dynamical dual --- \emph{bidynamical} Poisson groupoids.

Let $G$ be a connected, simply connected Lie group, $L\subset G$ a Lie
subgroup, and let $\l=\Lie(L)$ and $\g=\Lie(G)$ be their respective Lie
algebras; we denote by $i:\l\to\g$ the corresponding inclusion. For an
$\Ad^*_L$-equivariant subset $U\subset\l^*$ we consider the trivial Lie
groupoid $\GG=U\times G\times U$ with base $U$, with the product given
by $(p,x,q)(q,y,r)=(p,xy,r)$.  Let $r\colon U\to\bigwedge^2\g$ be a
differentiable map (we identify $\bigwedge^2\g$ with the skew-symmetric
maps from $\g^*$ to $\g$). In~\cite{EV98}, extending
\mbox{Drinfel\cprime{}d}'s classical work, P.~Etingof and A.~Varchenko,
showed that the following bracket on $C^\infty(\GG)$:
\begin{align*}
  \{f,g\}_{(p,x,q)}=&\,\<p,[\delta f,\delta g]_\l\>
  -\<q,[\delta'f,\delta'g]_\l\>
  -\<Dg,i\delta f\>-\<D'g,i\delta'f\>\\*
  &\quad
  +\<Df,i\delta g\>+\<D'f,i\delta'g\>
  -\<Df,r_pDg\>+\<D'f,r_qD'g\>
\end{align*}
is a Poisson bracket if and only if $r$ is a classical dynamical
$r$-matrix, in which case $\GG$ turns out to be a Poisson groupoid ---
which they call \emph{dynamical} (see Section~\ref{sc:Dyn at 0} for the
notations). We recall that a classical dynamical $r$-matrix is an
$\l$-equivariant solution of the classical dynamical Yang--Baxter
equation:
\begin{equation*}
  \Cycl_{(\xi,\eta,\zeta)}\Bigl(\<\zeta,\d_pr(i^*\xi)\eta\>
  -\<\zeta,[r_p\xi,r_p\eta]\>\Bigr)=
  \<\xi\otimes\eta\otimes\zeta,\varphi\>,
\end{equation*}
where $\varphi$ is any element in $\bigl(\bigwedge^3\g\bigr)^\g$. It
also appeared that the smallest Poisson submanifold of $\GG$ containing
the unit of $\GG$ is not the unit itself, but the image of the
hamiltonian unit $\LL=U\times L$ by the groupoid morphism
$I(p,h)=(\Ad^*_{h^{-1}}p,h,p)$. Explicit dynamical $r$-matrices were
given and classified when $\g$ is a complex semi-simple Lie algebra and
$\l$ a Cartan subalgebra. In~\cite{AM00}, A.~Alekseev and E.~Meinrenken
exhibited an analytic classical dynamical $r$-matrix in the case where
$\g=\l$ is a quadratic Lie algebra. In~\cite{ES01}, P.~Etingof and
O.~Schiffmann proved the existence of (formal) classical dynamical
$r$-matrices and gave a complete description of the moduli space of
classical dynamical $r$-matrices in the case where there exists a
reductive decomposition $\g=\l\oplus\m$ and for an element
$\varphi=\<\Omega,\Omega\>$ with $\Omega\in\bigl(\Sym^2\g\bigr)^\g$ such
that $\Omega\in\l\otimes\l\,\oplus\,\m\otimes\m$.

In~\cite{LP05}, L.C.~Li and S.~Parmentier gave the form of all Poisson
groupoid structures on the trivial Lie groupoid $\GG=U\times G\times U$
which admit an inclusion of the hamiltonian unit $\LL$. But in the
present paper, we only consider the following special form of Poisson
brackets on $\GG$ which correspond to the inclusion
$I(p,h)=(\Ad^*_{h^{-1}}p,h,p)$:
\begin{align*}
  \{f,g\}_{(p,x,q)}=&\,\<p,[\delta f,\delta g]_\l\>
  -\<q,[\delta'f,\delta'g]_\l\>
  -\<Dg,i\delta f\>-\<D'g,i\delta'f\>\\*
  &\quad+\<Df,i\delta g\>+\<D'f,i\delta'g\>
  -\<Df,l_pDg\>+\<Df,\pi_xDg\>+\<D'f,l_qD'g\>,
\end{align*}
where $\pi\colon G\to\bigwedge^2\g$ is a Lie group $1$-cocycle, and
$l:U\to\bigwedge^2\g$. It turns out that Jacobi's identity is equivalent
to the following two conditions:
\begin{itemize}
\item There exists an element $\varphi\in\bigwedge^3\g$ such that for
  all $\xi,\,\eta,\,\zeta\in\g^*$ the following identity holds:
  \begin{equation*}
    \<\xi\otimes\eta\otimes\zeta,\ad_x^{(3)}\varphi\>=
    \Cycl_{(\xi,\eta,\zeta)}
    \<\xi,\varpi_{\varpi_x\eta}\zeta\>,
  \end{equation*}
  and for all $p\in U$ and $\xi,\,\eta,\,\zeta\in\g^*$ the following
  identity holds:
  \begin{equation*}
    \Cycl_{(\xi,\eta,\zeta)}\Big(\<\zeta,\d_pl(i^*\xi)\eta\>
    -\<\zeta,[l_p\xi,l_p\eta]\>
    -\<\zeta,\varpi_{l_p\xi}\eta\>\Big)=
    \<\xi\otimes\eta\otimes\zeta,\varphi\>;
  \end{equation*}
\item for all $p\in U$, $z\in\l$, and $\xi\in\g^*$ the following
  identity holds:
  \begin{equation*}
    \d_pl(\ad^*_zp)\xi+\varpi_{iz}\xi+\ad_{iz}l_p\xi+l_p\ad^*_{iz}\xi=0.
  \end{equation*}
\end{itemize}
Such a groupoid $\GG$ is called \emph{dynamical}, and a map $l$
satisfying the two previous conditions will be referred to as a
\emph{classical dynamical $\ell$-matrix}.

It is shown in~\cite{PP05} that the previous conditions have solutions
$l\colon U\to\bigwedge^2\g$ only if the quadruple
$\G=(\g,[\,,\,],\varpi,\varphi)$ is a Lie quasi-bialgebra. From the form
of the dual of the Lie algebroid, it is observed that a necessary
condition for the groupoid $\GG$ to have a dynamical dual is that the
Lie algebra $\g$ admits a reductive decomposition $\g=\l\oplus\m$. Under
this assumption and additional natural (but restrictive) compatibility
conditions between $\G$ and the reductive decomposition, all formal
solutions $l$ of the above conditions are given, via an explicit and
analytic representative, and the action of a (formal) gauge group. The
duality for the groupoid $\GG$ is also explicitely described.

The goal of the present paper is to solve the problem when no
compatibility condition between the reductive decomposition of $\g$ and
the Lie quasi-bialgebra is assumed\message{description technique}. In
short Section~\ref{sc:Dyn at 0} we recall some notations and set the
problem. Section~\ref{sc:canonical l-matrices} is divided into three
subsections: first we show that every classical dynamical $\ell$-matrix
is gauge equivalent to one satisfying $l_psp=0$, and that there is at
most one formal classical dynamical $\ell$-matrix satisfying this
condition, which we call \emph{canonical}. Second we find an explicit
formula for the canonical $\ell$-matrix, which shows that it is analytic
(Theorem~\ref{th:analyticity}). Then, we define \emph{bidynamical}
objects and morphisms on both the Lie quasi-bialgebra level and groupoid
level, and show that there is a functorial correspondance between these
bidynamical objects. In Section~\ref{sec:TandD} we give an explicit
trivialization isomorphism which enables us to describe explicitly the
Poisson groupoid dual to $\GG$. The duality for groupoids induces a
duality for the class of bidynamical Lie quasi-bialgebras, which is also
described. We obtain the following characterization which was announced
in~\cite{PP05}: a dynamical Poisson groupoid (with $0\in U$) is
bidynamical if and only if $\g$ admits a reductive decomposition
$\g=\l\oplus\m$. In Section~\ref{sc:Link duality} we give a link between
the two canonical dynamical $\ell$-matrices associated to dynamical
Poisson groupoids in duality, which shows that both $\ell$-matrices have
the same domain of analyticity.

The problem will be adapted to the case where $0$ does not belong to $U$
in a forthcoming publication~\cite{P?}.

\begin{ack}
The author would like to thank S.\ Parmentier for useful discussions and
comments on this paper.
\end{ack}

\section{Dynamical Poisson groupoids at $0$}
\label{sc:Dyn at 0}
Let $i\colon\l\to\g$ be an inclusion of the Lie algebra $\l$ into the Lie
algebra $\g$. Let $U$ be an $\Ad^*_L$-invariant, contractible open
subset in $\l^*$ containing $0$ and $\G=(\g,[\,,\,],\varpi,\varphi)$ a
Lie quasi-bialgebra (see~\cite{PP05} for information on Lie
quasi-bialgebras where we use the same conventions and notations as in
the present paper).

By definition, a classical dynamical $\ell$-matrix on $U$ associated
with $\G$ is a map $l\colon U\to\bigwedge^2\g$ which is a solution of
the following two equations (we identify $\bigwedge^2\g$ with the
skew-symmetric maps from $\g^*$ to $\g$):
\begin{align}
  \label{eq:CDYBEg}
  \Cycl_{(\xi,\eta,\zeta)}\Big(\<\zeta,\d_pl(i^*\xi)\eta\>
  -\<\zeta,[l_p\xi,l_p\eta]\>
  -\<\zeta,\varpi_{l_p\xi}\eta\>\Big)
  &=\<\xi\otimes\eta\otimes\zeta,\varphi\>\\*
  \label{eq:l-eqg}
  \d_pl(\ad^*_zp)+\varpi_z+\ad_zl_p+l_p\ad^*_z&=0.
\end{align}
The set of classical dynamical $\ell$-matrices on $U$ associated with
$\G$ is denoted by $\Dynl(U,\G)$. Classical dynamical $r$-matrices are
just classical dynamical $\ell$-matrices with $\varpi=0$. Let $\DD$ be
the formal neighborhood of $0$ in $\l^*$. We also consider classical
dynamical $\ell$-matrices which are formal around $0\in\l^*$, the set of
which is denoted by $\Dynl(\DD,\G)$.

For all $\t\in\bigwedge^2\g$, we denote by $\G^\t$ the twist of the Lie
quasi-bialgebra $\G$ via $\t$. The following result is proved
in~\cite{PP05}:
\begin{prop}\label{pr:Dynltwist}
  For all $\t\in\bigwedge^2\g$,
  \begin{equation*}
    \Dynl(U,\G^\t)=\Dynl(U,\G)-\t.
  \end{equation*}
\end{prop}
This proposition allows us to be only concerned with classical dynamical
$\ell$-matrices which vanish at $0$, which form a set denoted by
$\Dynl_0(U,\G)$.

If we want $\Dynl_0(U,\G)$ (or $\Dynl_0(\DD,\G)$) to be non empty, then
we must have:
\begin{equation}
  \varpi_\l=0\qquad\text{and}\qquad
  \varphi\equiv0\mod\l
\end{equation}
(evaluate equations~\eqref{eq:CDYBEg} and~\eqref{eq:l-eqg} at $0$).

Classical dynamical $\ell$-matrices are related to Poisson groupoids in
the following way: Let $G$ be a connected Lie group with Lie algebra
$\g$. For any point $x\in G$ and any function $f\in C^\infty(G)$, we
denote by $D_xf\in\g^*$ and $D'_xf\in\g^*$ the right and left
derivatives at $x$:
\begin{align}
  \<D_xf,u\>&=\left.\frac\d{\d t}\right\vert_{t=0}f(\e^{tu}x)\\*
  \label{eq:D'_xf}
  \<D'_xf,u\>&=\left.\frac\d{\d t}\right\vert_{t=0}f(x\e^{tu})
\end{align}
for all $u\in\g$. Let $L$ be a connected Lie subgroup of $G$ with Lie
algebra $\l$, and $U$ an $\Ad^*_L$-invariant open subset in $\l^*$
containing $0$. We will denote the inclusion by
$i\colon\l\to\g$. Consider the trivial Lie groupoid $\GG=U\times G\times
U$ with multiplication:
\begin{equation}
  (p,x,q)(q,y,r)=(p,xy,r)
\end{equation}
We say that a multiplicative Poisson bracket on $\GG$ is
\emph{dynamical} if it is of the form:
\begin{equation}\label{eq:Dynamical_bracket}
\begin{aligned}
  \{f,g\}_{(p,x,q)}=&\,\<p,[\delta f,\delta g]_\l\>
  -\<q,[\delta'f,\delta'g]_\l\>\\*
  &\quad-\<Dg,i\delta f\>-\<D'g,i\delta'f\>\\*
  &\quad+\<Df,i\delta g\>+\<D'f,i\delta'g\>\\*
  &\quad-\<Df,l_pDg\>+\<Df,\pi_xDg\>+\<D'f,l_qD'g\>
\end{aligned}
\end{equation}
where $l\colon U\to\bigwedge^2\g$ is a smooth map, and $\pi\colon
G\to\bigwedge^2\g$ is a group $1$-cocycle. In this equation, $\delta
f\in\l$ and $\delta'f\in\l$ denote the derivatives of $f$ with respect
to the first and second $U$ factors, $Df$ and $D'f$ denote the right and
left derivatives of $f$ with respect to the $G$ factor, and all
derivatives are evaluated at $(p,x,q)$. Denote by
$\varpi=\T_1\pi\colon\g\to\bigwedge^2\g$ the Lie algebra $1$-cocycle
associated with $\pi$. It can be shown (see~\cite{LP05}) that the
bracket~\eqref{eq:Dynamical_bracket} is Poisson (\ie{}satisfies Jacobi's
identity) if and only if $l\in\Dynl(U,\G)$ with
$\G=(\g,[\,,\,],\varpi,\varphi)$ for some $\varphi\in\bigwedge^3\g$ such
that $\G$ is a Lie quasi-bialgebra.

There is a notion of duality for Poisson groupoids which extends that of
Poisson--Lie groups (see~\cite{M99,W88} and also~\cite{LP05,PP05} for
our more concrete case). It was already observed in~\cite{PP05} (see
also~\cite{LP05}) that the Lie algebra of the vertex group $G^\star_0$
is (isomorphic to) the lagrangian Lie subalgebra
$\g^\star_0=\l\oplus\l^\perp$ of the double $\dlie$ of the Lie
quasi-bialgebra $\G$, which is a reductive decomposition over $\l$ (we
recall that a \emph{reductive decomposition over $\mathfrak{b}$} of a
Lie algebra $\mathfrak{a}$ is a vector space decomposition
$\mathfrak{a}=\mathfrak{b}\oplus\mathfrak{c}$ such that $\mathfrak{b}$
is a Lie subalgebra of $\mathfrak{a}$ and
$[\mathfrak{b},\mathfrak{c}]\subset\mathfrak{c}$). Thus, a necessary
condition for the dual to be (a covering of) a dynamical Poisson
groupoid, is that $\g$ admits a reductive decomposition $\g=\l\oplus\m$.
In this case, it was shown in~\cite{PP05}, under additional natural but
restrictive compatibility conditions on $\varphi$ and $\varpi$, that the
dual of $\GG$ is still (a covering of) a dynamical Poisson groupoid.
This result is proved in this paper, without the additional
compatibility assumptions on $\varphi$ and $\varpi$. It is shown
in~\cite{LP05}, using a theorem of Mackenzie~\cite{M87} that the vertex
algebras $\g^\star_p$ defined as:
\begin{equation}\label{eq:dual_of_g}
  \g_p^\star=\ensemble{i(z)+\xi\in i(\l)\oplus\g^*}
  {i^*\xi=\ad^*_zp}\subset\g\oplus\g^*
\end{equation}
with the Lie bracket:
\begin{equation}\label{eq:bracket_on_gstar}
  \begin{aligned}\relax
    [i(z)+\xi,i(z')+\xi']^\star_p&=\bigl(i([z,z'])+\varpi_{i(z)}\xi'
    +\ad_{i(z)}l_p\xi' +l_p\ad_{i(z)}^*\xi'\\*
    &\quad-\varpi_{i(z')}\xi-\ad_{i(z')}l_p
    \xi-l_p\ad_{i(z')}^*\xi\\*
    &\quad+[l_p\xi,l_p\xi']+l_p\ad^*_{l_p\xi}\xi'
    -l_p\ad^*_{l_p\xi'}\xi\\*
    &\quad+\varpi_{l_p\xi}\xi'-\varpi_{l_p\xi'}\xi
    -\<\xi,\varpi_{l_p\bullet}\xi'\>
    +\<\xi\otimes \xi'\otimes1,\varphi\>,\\*
    &\quad-\ad^*_{i(z)}\xi'+\ad^*_{i(z')}\xi-\<\xi,\varpi_\bullet\xi'\>
    -\ad^*_{l_p\xi}\xi'+\ad^*_{l_p\xi'}\xi\bigr)
  \end{aligned}
\end{equation}
are all isomorphic when $p$ ranges over $U$. Isomorphisms between the
$\g^\star_p$'s are not provided by Mackenzie's theorem, and are not
canonical. In the present paper, we construct such an isomorphism as
part of the trivialization. It was observed in~\cite{PP05} that the Lie
algebra $\g^\star_p$ is a lagrangian Lie subalgebra in $\dlie^p$, the
canonical double of the twisted Lie quasi-bialgebra $\G^{l_p}$.

It is natural to consider the following definition:
\begin{defn}
  A dynamical Poisson groupoid $\GG$ over $U$ is said to be
  \emph{bidynamical} if its dual is (a covering of) a dynamical Poisson
  groupoid.
\end{defn}
\begin{rem}
  It is shown in~\cite{PP05} that if the Lie quasi-bialgebra $\G$ is
  compatible with the reductive decomposition $\g=\l\oplus\m$ in the
  sense of Example~\ref{ex:lmatr}, the associated dynamical Poisson
  groupoid is bidynamical.
\end{rem}

This paper is devoted to the study of bidynamical Poisson groupoids (at
$0$). Thus, from now on, we fix a reductive decomposition
$\g=\l\oplus\m$ of a Lie algebra $\g$, and a Lie quasi-bialgebra
$\G=(\g,[\,,\,],\varpi,\varphi)$ such that $\varpi_\l=0$ and
$\varphi\equiv0\mod\l$. We also denote by $\p_\l$ (resp.~$\p_\m$) the
projection on $\l$ (resp.~$\m$) along $\m$ (resp.~$\l$), and by
$s\colon\l^*\to\g^*$ the adjoint of $\p_\l$. Notice that $s$ is
$\l$-equivariant and that its image is $\m^\perp$.

\section{Canonical dynamical $\ell$-matrices}
\label{sc:canonical l-matrices}
For a vector space $E$ and a formal map $f\colon\DD\to E$, we denote
either by $[f]_k$ or $f^k$ its homogeneous term of degree $k$ (the
notation $f^k$ is not to be confused with the $k$-th power of $f$).

For a dynamical $\ell$-matrix $l\in\Dynl(\DD,\G)$ and an
$\l$-equivariant map $\sigma\in\Map(\DD,G)^\l$, we denote by $l^\sigma$
the gauge transformation of the map $l$ by $\sigma$ (see~\cite{PP05}):
\begin{equation}\label{eq:actionsurdynl}
  l^\sigma_p=\Ad_{\sigma_p}l_p\Ad^*_{\sigma_p}+\theta^\sigma_p
  +\pi_{\sigma_p},
\end{equation}
where $\theta^\sigma$ is defined as:
\begin{equation}
  \theta^\sigma_p=r_{\sigma_p^{-1}}(\T_p\sigma)i^*\Ad_{\sigma_p}^*
  -(\T_p\sigma)^*r^*_{\sigma_p^{-1}},
\end{equation}
and where $\pi\colon G\to\bigwedge^2\g$ is the Lie group cocycle
integrating $\varpi$. This action is a left action:
\begin{equation}
  (l^\sigma)^{\sigma'}=l^{\sigma'\sigma}.
\end{equation}

\subsection{Canonical dynamical $\ell$-matrices}
The following proposition shows that there always exists some
distinguished representative in each (formal) gauge orbit of dynamical
$\ell$-matrices.

We denote by $\Map(\DD,G)$ the group of formal maps from $\DD$ to $G$
with pointwise multiplication, by $\Map_0(\DD,G)$ the group of formal
maps $\sigma$ such that $\sigma_0=1$ and by $\Map^{(2)}_0(\DD,G)$ the
group of formal maps $\sigma$ such that $\sigma_0=1$ and $\T_0\sigma=0$.
\begin{prop}\label{pr:lpsp=0}
  For all $l\in\Dynl_0(\DD,\G)$, there exists a gauge transformation
  associated with some $\sigma\in\Map^{(2)}_0(\DD,G)^\l$ such that
  $l^\sigma_psp=0$. 
\end{prop}
\begin{proof}
  By induction: Let $k\geq1$ and assume that $[l_p]_{\leq k-1}sp=0$. 
  Now, if $\Sigma\colon\DD\to\g$ is an $\l$-equivariant homogeneous map
  of degree $k+1$, then setting $\sigma=\e^\Sigma$ yields
  $[l^\sigma]_k=[l]_k+\d\Sigma i^*-(\d\Sigma)^*$, so that
  $[l_p^\sigma]_{\leq k}sp=[l_p]_ksp+\d_p\Sigma(p)-(\d_p\Sigma)^*(sp)$. 
  Now, define $\Sigma$ as:
  \begin{equation}\label{eq:Sigma de lpsp=0}
    \Sigma_p=-\frac1{k+1}\p_\m[l_p]_ksp-\frac1{k+2}\p_\l[l_p]_ksp.
  \end{equation}
  Notice that $\Sigma$ is an $\l$-equivariant homogeneous map of degree
  $k+1$. We have 
  \begin{align*}
    \d_p\Sigma(\alpha)&=
    -\left(\frac1{k+1}\p_\m+\frac1{k+2}\p_\l\right)
    \bigl([l_p]_ks\alpha+\d_p[l]_k(\alpha)sp\bigr)
    \intertext{thus}
    \d_p\Sigma(p)&=-\p_\m[l_p]_ksp-\frac{k+1}{k+2}[l_p]_ksp
    \intertext{and}
    (\d_p\Sigma)^*sp&=\frac1{k+2}\p_\l[l_p]_ksp.
  \end{align*}
  Hence, $\d_p\Sigma(p)-(\d_p\Sigma)^*sp=-[l_p]_ksp$ and
  $[l^\sigma_p]_{\leq k}sp=0$. The proof follows by induction. 
\end{proof}

We will make use of the following notations: for all $\xi\in\g^*$, we
define:
\begin{equation}
  \xi'_p=\p_{\g^*}\ad_{sp}l_p\xi,\qquad
  \xi''_p=\p_{\g^*}\ad_{sp}\xi,\qquad
  \widetilde{\xi}_p=\xi'_p+\xi''_p.
\end{equation}
The following proposition shows that such a representative is
necessarily unique.
\begin{prop}\label{pr:uniqueness}
  There exists at most one (formal) dynamical $\ell$-matrix $l$
  satisfying $l_psp=0$ and $l_0=0$. It is the unique (formal) solution
  $l$ satisfying $l_0=0$ of the following differential equation:
  \begin{equation}\label{eq:diffeq pour lcan}
    \d_pl(p)\zeta=\ad_{sp}(l_p\zeta+\zeta)-
    (l_p\widetilde\zeta_p+\widetilde\zeta_p)-\p_\l l_p\zeta-l_psi^*\zeta
  \end{equation}
  for all $\zeta\in\g^*$.
\end{prop}
\begin{proof}
Assume that $l\colon\DD\to\bigwedge^2\g$ is a formal map satisfying
$l_psp=0$ and the generalized Yang--Baxter equation~\eqref{eq:CDYBEg}.

For all $\xi,\,\eta\in\l^\perp$, we obtain from
equation~\eqref{eq:CDYBEg}:
\begin{align}
  \nonumber
  \<\eta,\d_pl(p)\xi\>&=
  (sp,[l_p\xi,l_p\eta+\eta])_\dlie+(\xi,[l_p\eta,sp])_\dlie
  +(\eta,\ad_{sp}\xi)_\dlie\\*
  \label{eq:15264856}
  &=(\eta,-l_p\widetilde{\xi}_p+\ad_{sp}l_p\xi+\ad_{sp}\xi)_\dlie.
\end{align}
If we define $L\colon\DD\to\L(\g^*,\g)$ as:
\begin{equation}
  L_p\zeta=\p_\g\ad_{sp}\zeta+\ad_{sp}l_p\zeta-l_p\widetilde{\zeta}_p
\end{equation}
for all $\zeta\in\g^*$, then equation~\eqref{eq:15264856} reads:
\begin{equation}
  \p_\m[l_p]_k\xi=\frac1k\p_\m[L_p]_k\xi
\end{equation}
for all $\xi\in\l^\perp$ and $k\geq1$. 

For all $\xi\in\l^\perp$ and $\alpha\in\l^*$, we obtain from
equation~\eqref{eq:CDYBEg}:
\begin{equation}\label{eq:152648562}
  \<\xi,\d_pl(\alpha)sp-\d_pl(p)s\alpha\>=\<\xi,L_ps\alpha\>.
\end{equation}
Since $\d_pl(\alpha)sp=-l_ps\alpha$ for all $\alpha\in\l^*$,
equation~\eqref{eq:152648562} reads:
\begin{equation}
  \p_\m[l_p]_ks\alpha=\frac1{k+1}\p_\m[L_p]_ks\alpha
\end{equation}
for all $\alpha\in\l^*$ and $k\geq1$. 
Similarly, one obtains:
\begin{align}
  \p_\l[l_p]_k\xi&=\frac1{k+1}\p_\l[L_p]_k\xi\\*
  \p_\l[l_p]_ks\alpha&=\frac1{k+2}\p_\l[L_p]_ks\alpha
\end{align}
for all $\xi\in\l^\perp$, $\alpha\in\l^*$ and $k\geq1$. Thus, $[l_p]_k$
is uniquely determined by the $[l_p]_j$'s, for $j\le k-1$. It is easily
shown that $l$ satisfies equation~\eqref{eq:diffeq pour lcan}.
\end{proof}

In particular, one obtains the following corollary:
\begin{cor}
  The reduced moduli space $\Dynl_0(\DD,\G)/\Map_0^{(2)}(\DD,G)^\l$
  consists of at most one point. 
\end{cor}
\begin{proof}
  Let $l$ and $l'$ in $\Dynl_0(\DD,\G)$. Then, from
  Proposition~\ref{pr:lpsp=0} there are two maps $\sigma$ and $\sigma'$
  in $\Map_0^{(2)}(\DD,G)^\l$ such that
  $l^\sigma_psp=(l')^{\sigma'}_psp=0$. Thus, $l^\sigma$ and
  $(l')^{\sigma'}$ are two solutions of equation~\eqref{eq:diffeq pour
    lcan} which vanish at zero. From the uniqueness of such solutions, we
  must have $l^\sigma=(l')^{\sigma'}$. 
\end{proof}

Conversely, we show that a solution of~\eqref{eq:diffeq pour lcan} is a
(formal) dynamical $\ell$-matrix:
\begin{thm}\label{th:l-matrices generales}
  The (formal) solution $l$ of~\eqref{eq:diffeq pour lcan} with initial
  condition $l_0=0$ is the unique (formal) dynamical $\ell$-matrix
  satisfying $l_0=0$ and $l_psp=0$. 
\end{thm}
Before proving Theorem~\ref{th:l-matrices generales}, we introduce some
notations and state three lemmas:

\begin{lem}\label{lm:solution de diffeq}
  The (formal) solution $l$ of equation~\eqref{eq:diffeq pour lcan} such
  that $l_0=0$ takes its values in $\bigwedge^2\g$, and satisfies
  equation~\eqref{eq:l-eqg}.
\end{lem}
\begin{proof}
  It is a direct consequence of the skew-symmetry of the map $\ad_{sp}$,
  and of the $\l$-equivariance of the maps $s$ and $\p_\l$.
\end{proof}

\begin{lem}\label{lm:xitilde}
  For all $\xi\in\g^*$, the element $\widetilde\xi_p+\p_\l l_p\xi$ lies
  in $\g^\star_p$. 
\end{lem}
\begin{proof}
  For all $z\in\l$,
  $\<\widetilde\xi_p,z\>=\bigl(\ad_{sp}(l_p\xi+\xi),z\bigr)_\dlie
  =\bigl(\ad_{sp}l_p\xi,z\bigr)_\dlie =\bigl(\ad_{sp}\p_\l
  l_p\xi,z\bigr)_\dlie$, since $\varpi_\l=0$. Thus,
  $i^*\widetilde\xi_p=\ad^*_{\p_\l l_p\xi}p$, and Lemma~\ref{lm:xitilde}
  is proved. 
\end{proof}

We denote by $\CDYB$ the \emph{(generalized) classical dynamical
  Yang--Baxter operator}: For any map $l\in\Map(\DD,\bigwedge^2\g)$,
$\CDYB(l)\colon\DD\to\bigwedge^3\g$ is given by
\begin{equation}
  \<\xi\otimes\eta\otimes\zeta,\CDYB_p(l)\>=
  \Cycl_{(\xi,\eta,\zeta)}\big(\zeta,\d_pl(i^*\xi)\eta-
  [l_p\xi,l_p\eta+\eta]\big)_\dlie
\end{equation}
for all $\xi,\,\eta,\,\zeta\in\g^*$, where, as usual, $(\,,\,)_\dlie$ is
the canonical bilinear form on the double $\dlie$ of $\G$. Notice that
equation~\eqref{eq:CDYBEg} is equivalent to $\CDYB(l)=\varphi$.

\begin{lem}\label{lm:CDY be-bop}
  For all $\xi,\,\eta,\,\zeta\in\g^*$, we have
  \begin{equation}
    \d_p\<\xi\otimes\eta\otimes\zeta,\CDYB_\bullet(l)\>(p)=
    -\Cycl_{(\xi,\eta,\zeta)}
    \<\xi\otimes(si^*\eta+\widetilde\eta_p)\otimes\zeta,
    \CDYB_p(l)-\varphi\>.
  \end{equation}
\end{lem}
\begin{proof}
  For all $\xi,\,\eta,\,\zeta\in\g^*$, using equations~\eqref{eq:diffeq
    pour lcan} and~\eqref{eq:l-eqg}, as well as:
  \begin{align}
    \Cycl_{(\xi,\eta,\zeta)}\<\xi'_p,[l_p\eta,l_p\zeta]\>&=0\\*
    \Cycl_{(\xi,\eta,\zeta)}\<\xi''_p,[\eta,\zeta]\>&=0
  \end{align}
  which are consequences of Jacobi's identity on $\dlie$, one obtains:
  \begin{align}
    \nonumber
    \begin{split}
      \d_p&\Cycl_{(\xi,\eta,\zeta)}
      \bigr(\xi,[l\eta,l\zeta+\zeta]\big)_\dlie(p)
      =\Cycl_{(\xi,\eta,\zeta)}\big(\xi,[\d_pl(p)\eta,l_p\zeta+\zeta]
      +[l_p\eta,\d_pl(p)\zeta]\bigl)_\dlie\\*
      &=\Cycl_{(\xi,\eta,\zeta)}\big(\xi,
      [\ad_{sp}l_p\eta,l_p\zeta+\zeta]+[\ad_{sp}\eta,l_p\zeta]
      -[l_p\widetilde\eta_p+\widetilde\eta_p,l_p\zeta+\zeta]
      -[\p_\l l_p\eta,\zeta]\\*
      &\qquad -[\p_\l l_p\eta,l_p\zeta] -[l_psi^*\eta,l_p\zeta+\zeta]
      +[l_p\eta,\ad_{sp}(l_p\zeta+\zeta)
      -(l_p\widetilde\zeta_p+\widetilde\zeta_p) -\p_\l
      l_p\zeta-l_psi^*\zeta]
      \big)_\dlie\\*
      &=\Cycl_{(\xi,\eta,\zeta)}\big(\xi,
      \ad_{sp}[l_p\eta,l_p\zeta+\zeta]+[\ad_{sp}\eta,l_p\zeta]
      -[l_p\widetilde\eta_p,l_p\zeta+\zeta]
      -[\widetilde\eta_p,l_p\zeta+\zeta]\\*
      &\qquad
      +\d_pl(i^*\widetilde\eta_p)\zeta-[l_psi^*\eta,l_p\zeta+\zeta]
      -[l_p\eta,l_p\widetilde\zeta_p+\widetilde\zeta_p]
      -[l_p\eta,l_psi^*\zeta]\big)_\dlie\\*
      &=\Cycl_{(\xi,\eta,\zeta)}
      -\big(\widetilde\eta_p,[l_p\zeta,l_p\xi]\big)_\dlie
      -\big(\xi,[l_p\widetilde\eta_p,l_p\zeta+\zeta]\big)_\dlie
      -\big(\widetilde\eta_p,[l_p\zeta,\xi]+[\zeta,\xi]\big)_\dlie\\*
      &\qquad +\big(\xi,\d_pl(i^*\widetilde\eta_p)\zeta
      -[l_psi^*\eta,l_p\zeta+\zeta]\big)_\dlie
      -\big(\zeta,[l_p\xi,l_p\widetilde\eta_p+\widetilde\eta_p]
      -[l_p\xi,l_psi^*\eta]\big)_\dlie
    \end{split}\\*[.5em]
    \label{eq:gre56g114e56gvb}
    \begin{split}
      &=\Cycl_{(\xi,\eta,\zeta)}
      \<\xi\otimes\widetilde\eta_p\otimes\zeta,\CDYB_p(l)-\varphi\>
      -\<\widetilde\eta_p,\d_pl(i^*\zeta)\xi\>-
      \<\zeta,\d_pl(i^*\xi)\widetilde\eta_p\>\\*
      &\qquad
      -\<\zeta,[l_p\xi,l_psi^*\eta]\>
      -\big(\xi,[l_psi^*\eta,l_p\zeta+\zeta]\big)_\dlie.
    \end{split}
  \end{align}

  Now,
  \begin{align*}
    \<\xi,l_p\zeta\>&=\sum_{k\geq0}[\<\xi,l_p\zeta\>]_k
    =\<\xi,l_0\zeta\>+\sum_{k\geq1}\frac1k[\<\xi,d_pl(p)\zeta\>]_k\\*
    &=\<\xi,l_0\zeta\>+\sum_{k\geq1}\frac1k\left[
      \big(\xi,\ad_{sp}(l_p\zeta+\zeta)
      -l_p\widetilde\zeta_p-\p_\l l_p\zeta-l_psi^*\zeta\big)_\dlie
      \right]_k
  \end{align*}
  thus,
  \begin{align*}
    \Cycl_{(\xi,\eta,\zeta)}\<\xi,\d_pl(i^*\eta)\zeta\>&=
    \Cycl_{(\xi,\eta,\zeta)}\sum_{k\geq1}\frac1k\Big[
    \<\xi,\ad_{si^*\eta}(l_p\zeta+\zeta)+\ad_{sp}\d_pl(i^*\eta)\zeta
    -\d_pl(i^*\eta)\widetilde\zeta_p-l_p\p_{g^*}\ad_{si^*\eta}l_p\zeta\\*
    &\qquad\qquad
    -l_p\p_{\g^*}\ad_{sp}\d_pl(i^*\eta)\zeta-l_p\p_{\g^*}\ad_{i^*\eta}\zeta
    -\p_\l\d_pl(i^*\eta)\zeta-\d_pl(i^*\eta)si^*\zeta\>
    \Big]_{k-1}\\*
    &=\Cycl_{(\xi,\eta,\zeta)}\sum_{k\geq1}\frac1k\Big[
    \big(si^*\eta,[l_p\zeta,l_p\xi+\xi]+[\zeta,\xi]\big)_\dlie
    -\<\widetilde\xi_p,\d_pl(i^*\eta)\zeta\>
    -\<\xi,\d_pl(i^*\eta)\widetilde\zeta_p\>\\*
    &\qquad\qquad
    +\big(\zeta,[l_p\xi,si^*\eta]\big)_\dlie
    -\<si^*\xi,\d_pl(i^*\eta)\zeta\>-\<\xi,\d_pl(i^*\eta)si^*\zeta\>
    \Big]_{k-1}
  \end{align*}
  and we obtain
  \begin{align}
    \nonumber
    \Cycl_{(\xi,\eta,\zeta)}\d_p\<\xi,\d_\bullet l(i^*\eta)\zeta\>(p)&=
    \Cycl_{(\xi,\eta,\zeta)}
    -\<\xi,\d_pl(i^*\eta)\zeta\>
    +\bigl(si^*\eta,[l_p\zeta,l_p\xi+\xi]\bigr)_\dlie
    -\<\widetilde\eta_p,\d_pl(i^*\zeta)\xi\>\\*
    \nonumber
    &\qquad
    -\<\zeta,\d_pl(i^*\xi)\widetilde\eta_p\>
    +\bigl(\zeta,[l_p\xi,si^*\eta]\bigr)_\dlie
    -\<si^*\eta,\d_pl(i^*\zeta)\xi\>\\*
    \nonumber
    &\qquad
    -\<\zeta,\d_pl(i^*\xi)si^*\eta\>
    +\bigl(\xi,[si^*\eta,\zeta]\bigr)_\dlie\\*
    \begin{split}
    \label{eq:res4rg61sebgr87s6e}
    &=\Cycl_{(\xi,\eta,\zeta)}
    -\<\xi\otimes i^*\eta\otimes\zeta,\CDYB_p(l)-\varphi\>
    -\<\zeta,[l_p\xi,l_psi^*\eta]\>\\*
    &\qquad
    -\bigl(\xi,[l_psi^*\eta,l_p\zeta+\zeta]\bigr)_\dlie
    -\<\widetilde\eta_p,\d_pl(i^*\zeta)\xi\>
    -\<\zeta,\d_pl(i^*\xi)\widetilde\eta_p\>.
    \end{split}
  \end{align}
  Assembling equations~\eqref{eq:gre56g114e56gvb}
  and~\eqref{eq:res4rg61sebgr87s6e} proves Lemma~\ref{lm:CDY be-bop}. 
\end{proof}

We can now prove Theorem~\ref{th:l-matrices generales}:
\begin{proof}[Proof of Theorem~\ref{th:l-matrices generales}]
  The map $l$ is $\l$-equivariant and takes its values in
  $\bigwedge^2\g$, by Lemma~\ref{lm:solution de diffeq}. We show that
  $l$ satisfies the generalized Yang--Baxter equation~\eqref{eq:CDYBEg}
  by induction:
  
  A computation shows that it is satified at order $0$,
  \ie~$\CDYB_0(l)=\varphi$. Indeed, if $\xi,\,\eta,\,\zeta\in\l^\perp$,
  then $\Cycl_{(\xi,\eta,\zeta)}\<\xi,\d_0l(i^*\eta)\zeta\>=0$ and
  \begin{equation}
    \Cycl_{(\xi,\eta,\zeta)}\bigl(\xi,[l_0\eta,l_0\zeta+\zeta]\bigr)_\dlie
    =0=\<\xi\otimes\eta\otimes\zeta,\varphi\>
  \end{equation}
  since $\varphi\equiv0\mod\l$. Thus,
  $\p_\m^{(3)}\CDYB_0(l)=\p_\m^{(3)}\varphi$. Now, from
  Lemma~\ref{lm:CDY be-bop}, if $\alpha,\,\beta,\,\gamma\in\l^*$, we
  obtain:
  \begin{align}
    0&=\<s\alpha\otimes\xi\otimes\eta,\CDYB_0(l)-\varphi\>\\*
    0&=2\<s\alpha\otimes s\beta\otimes\xi,\CDYB_0(l)-\varphi\>\\*
    0&=3\<s\alpha\otimes s\beta\otimes s\gamma,\CDYB_0(l)-\varphi\>.
  \end{align}
  Thus, $\CDYB_0(l)=\varphi$. 

  So let $k\geq0$ and assume that $[\CDYB_p(l)-\varphi]_{\leq k}=0$. 
  From Lemma~\ref{lm:CDY be-bop}, since $l$ satisfies
  equation~\eqref{eq:diffeq pour lcan}, we have:
  \begin{equation}\label{eq:ghs4e56ht4rse56h}
    (k+1)[\<\xi\otimes\eta\otimes\zeta,\CDYB_p(l)\>]_{k+1}=
    -\Cycl_{(\xi,\eta,\zeta)}[\<\xi\otimes(si^*\eta+\widetilde\eta_p)
    \otimes\zeta,\CDYB_p(l)-\varphi\>]_{k+1}
  \end{equation}
  for all $\xi,\,\eta,\,\zeta\in\g^*$. Now, since $\widetilde\eta_0=0$,
  equation~\eqref{eq:ghs4e56ht4rse56h} reads:
  \begin{equation}\label{eq:graou !!!} 
    (k+1)[\<\xi\otimes\eta\otimes\zeta,\CDYB_p(l)\>]_{k+1}=
    -\Cycl_{(\xi,\eta,\zeta)}[\<\xi\otimes si^*\eta\otimes\zeta,
    \CDYB_p(l)\>]_{k+1}
  \end{equation}
  for all $\xi,\,\eta,\,\zeta\in\g^*$. 
  \begin{enumerate}
  \item Let $\xi,\,\eta,\,\zeta\in\l^\perp$. Then
    equation~\eqref{eq:graou !!!} reads
    $[\<\xi\otimes\eta\otimes\zeta,\CDYB_p(l)\>]_{k+1}=0$;
  \item let $\xi,\,\eta\in\l^\perp$ and $\zeta\in\m^\perp$. Then
    equation~\eqref{eq:graou !!!} reads:
    \begin{equation*}
      (k+1)[\<\xi\otimes\eta\otimes\zeta,\CDYB_p(l)\>]_{k+1}=
      -[\<\xi\otimes\eta\otimes\zeta,\CDYB_p(l)\>]_{k+1}
    \end{equation*}
    thus $[\<\xi\otimes\eta\otimes\zeta,\CDYB_p(l)\>]_{k+1}=0$;
  \item let $\xi\in\l^\perp$ and $\eta,\,\zeta\in\m^\perp$. Then
    equation~\eqref{eq:graou !!!} reads:
    \begin{equation*}
      (k+1)[\<\xi\otimes\eta\otimes\zeta,\CDYB_p(l)\>]_{k+1}=
      -2[\<\xi\otimes\eta\otimes\zeta,\CDYB_p(l)\>]_{k+1}
    \end{equation*}
    thus $[\<\xi\otimes\eta\otimes\zeta,\CDYB_p(l)\>]_{k+1}=0$;
  \item let $\xi,\,\eta,\,\zeta\in\m^\perp$. Then
    equation~\eqref{eq:graou !!!} reads:
    \begin{equation*}
      (k+1)[\<\xi\otimes\eta\otimes\zeta,\CDYB_p(l)\>]_{k+1}=
      -3[\<\xi\otimes\eta\otimes\zeta,\CDYB_p(l)\>]_{k+1}
    \end{equation*}
    thus $[\<\xi\otimes\eta\otimes\zeta,\CDYB_p(l)\>]_{k+1}=0$. 
  \end{enumerate}
  Hence, the map $l$ satisfies equation~\eqref{eq:CDYBEg}. Thus
  $l\in\Dynl(\DD,\G)$. We have already seen that it is the unique
  dynamical $\ell$-matrix satisfying $l_0=0$ and $l_psp=0$. 
\end{proof}

\begin{defn}\label{df:Gen l-matrix}
  The (formal) $\ell$-matrix defined in Theorem~\ref{th:l-matrices
    generales} is called the \emph{canonical dynamical $\ell$-matrix
    associated with the Lie quasi-bialgebra $\G$ and the reductive
    decomposition $\g=\l\oplus\m$}, and is denoted by
  $\lcanG{\G,\l,\m}$, or simply $\lcan$ when no confusion is possible. 
\end{defn}

\begin{rem}
  Clearly, when the Lie quasi-bialgebra $\G$ is canonically compatible
  with the reductive decomposition $\g=\l\oplus\m$ of $\g$ (see
  Example~\ref{ex:lmatr} or~\cite{PP05}), then the canonical
  $\ell$-matrices of Definition~\ref{df:Gen l-matrix} and of~\cite{PP05}
  coincide in view of uniqueness of dynamical $\ell$-matrices satisfying
  $l_psp=0$ (see also Example~\ref{ex:lmatr}).
\end{rem}

\subsection{Canonical dynamical $\ell$-matrices are analytic}
We now show that a canonical dynamical $\ell$-matrix is in fact
analytic, and find an explicit formula.

For any $\t\in\bigwedge^2\g$ we set
\begin{equation*}
  \appli{\tau_\t}{\g\oplus\g^*}{\g\oplus\g^*}
  {x+\xi}{(x+\t\xi)+\xi.}
\end{equation*}

 We start with the following
proposition:
\begin{prop}\label{pr:pmadtauzero}
  Let $l=\lcanG{\G,\l,\m}\in\Dynl_0(\DD,\G)$. Then, for all
  $X_p\in\g^\star_p$ and for all $\alpha\in\l^*$, one has:
  \begin{align}
    \label{eq:pmadtauzero}
    \p_\m\Ad_{\e^{-sp}}\tau_{l_p}X_p&=0\\*
    \label{eq:pmdadtauzero}
    \p_\l\frac{\Ad_{\e^{-sp}}-1}{\ad_{sp}}\tau_{l_p}X_p&=-\p_\l X_p\\*
    \label{eq:pmadalpha}
    \p_\m\Ad_{\e^{-sp}}\tau_{l_p}s\alpha&=
    -\p_\m\frac{\Ad_{\e^{-sp}}-1}{\ad_{sp}}s\alpha\\*
    \label{eq:pldadtaualpha}
    \p_\l\frac{\Ad_{\e^{-sp}}-1}{\ad_{sp}}\tau_{l_p}s\alpha&=
    -\p_\l\frac{\Ad_{\e^{-sp}}-1+\ad_{sp}}{\ad_{sp}^2}s\alpha.
  \end{align}
\end{prop}
We recall that $X_p$ belongs to $\g^\star_p$ if and only if
$X_p=z_p+\xi_p\in\l\oplus\g^*$ and $i^*\xi_p=\ad^*_{z_p}p$.
\begin{proof}
  Since $l_0=0$, equations~\eqref{eq:pmadtauzero}
  and~\eqref{eq:pmdadtauzero} are clearly satisfied at order $0$. Now,
  let $X_p\in\g^\star_p$. Then, $X_p$ can be uniquely written
  $X_p=z+\ad_{sp}z+\xi$, where $z\in\l$ and $\xi\in\l^\perp$ (notice
  that $\ad_{sp}z=\ad^*_zsp\in\m^\perp$ for all $z\in\l$, $p\in\l^*$,
  since $\varpi_\l=0$). 

  Since $l$ satisfies equation~\eqref{eq:diffeq pour lcan}, one has:
  \begin{align*}
    \d_p(\Ad_{\e^{-s\bullet}}\tau_{l_\bullet}X_\bullet)(p)&=
    \Ad_{\e^{-sp}}\Bigl(-\ad_{sp}\tau_{l_p}X_p
    +\ad_{sp}\tau_{l_p}(X_p-z)
    -\tau_{l_p}\bigl(\widetilde\xi_p+\widetilde{(\ad_{sp}z)}_p\bigr)\\*
    &\qquad
    -\p_\l l_p(X_p-z)-l_p\ad_{sp}z+\tau_{l_p}\ad_{sp}z\Bigr)\\*
    &=-\Ad_{\e^{-sp}}
    \tau_{l_p}\bigl(\widetilde\xi_p+\widetilde{(\ad_{sp}z)}_p
    +\p_\l l_p(\xi+\ad_{sp}z)\bigr).
  \end{align*}
  Now, the element $\widetilde\xi_p+\widetilde{(\ad_{sp}z)}_p+\p_\l
  l_p(\xi+\ad_{sp}z)$ lies in $\g^\star_p$, by Lemma~\ref{lm:xitilde},
  and is of degree $\geq1$. Thus, if equation~\eqref{eq:pmadtauzero} is
  satisfied modulo terms of degree $\geq k$, then it will also be
  satisfied modulo terms of degree $\geq k+1$. 
  Equation~\eqref{eq:pmadtauzero} is thus proved. 
  
  Since $l$ satisfies equation~\eqref{eq:diffeq pour lcan}, one has:
  \begin{align*}
    \d_p\left(\frac{\Ad_{\e^{-s\bullet}}-1}{\ad_{s\bullet}}
      \tau_{l_\bullet}X_\bullet\right)(p)&=
    -\Ad_{\e^{-sp}}\tau_{l_p}X_p
    -\frac{\Ad_{\e^{-sp}}-1}{\ad_{sp}}\tau_{l_p}X_p\\*
    &\qquad
    +\frac{\Ad_{\e^{-sp}}-1}{\ad_{sp}}
    \Bigl(\ad_{sp}\tau_{l_p}X_p
    -\tau_{l_p}\bigl(\widetilde\xi_p+\widetilde{(\ad_{sp}z)}_p\bigr)
    -\p_\l l_p(X_p-z)\Bigr)\\*
    &=-\tau_{l_p}X_p-\frac{\Ad_{\e^{-sp}}-1}{\ad_{sp}}\tau_{l_p}\bigl(
    \widetilde\xi_p+(\widetilde{\ad_{sp}z})_p+\p_\l(\xi+\ad_{sp}z)+X_p
    \bigr).
  \end{align*}
  Thus, for $k\geq1$,
  \begin{align*}
    (k+1)\left[\p_\l\frac{\Ad_{\e^{-sp}}-1}{\ad_{sp}}\tau_{l_p}X_p\right]_k&=
    -[\p_\l l_p(\xi+\ad_{sp}z)]_k\\*
    &\qquad
    -\left[
      \p_\l\frac{\Ad_{\e^{-sp}}-1}{\ad_{sp}}\tau_{l_p}\bigl(
      \widetilde\xi_p+(\widetilde{\ad_{sp}z})_p+\p_\l(\xi+\ad_{sp}z)
      \bigr)
    \right]_k.
  \end{align*}
  Now, assume that equation~\eqref{eq:pmdadtauzero} is satisfied modulo
  terms of degree $\geq k$ (for some $k\geq1$). Then,
  \begin{equation*}
    (k+1)\left[\p_\l
      \frac{\Ad_{\e^{-sp}}-1}{\ad_{sp}}\tau_{l_p}X_p\right]_k=0
  \end{equation*}
  and equation~\eqref{eq:pmdadtauzero} is satisfied modulo terms of
  degree $\geq k+1$, and equation~\eqref{eq:pmdadtauzero} is proved. 
  
  Let $\alpha\in\l^*$. Similarly, one obtains:
  \begin{align*}
    (k+1)\left[
      \Ad_{\e^{-sp}}(l_ps\alpha+s\alpha)
      +\frac{\Ad_{\e^{-sp}}-1}{\ad_{sp}}s\alpha
    \right]_k&=-\left[\Ad_{\e^{-sp}}\tau_{l_p}(\widetilde{s\alpha}_p
    +\p_\l l_ps\alpha)\right]_k.
  \end{align*}
  Thus, by equation~\eqref{eq:pmadtauzero}:
  \begin{equation*}
    \p_\m\left(
      \Ad_{\e^{-sp}}(l_ps\alpha+s\alpha)
      +\frac{\Ad_{\e^{-sp}}-1}{\ad_{sp}}s\alpha
    \right)=0
  \end{equation*}
  and equation~\eqref{eq:pmadalpha} is proved. 

  Let $\alpha\in\l^*$. Similarly, one obtains:
  \begin{align*}
    (k+2)\left[
      \frac{\Ad_{\e^{-sp}}-1}{\ad_{sp}}\tau_{l_p}s\alpha
      +\frac{\Ad_{\e^{-sp}}-1+\ad_{sp}}{\ad_{sp}^2}s\alpha
    \right]_k&=-[\tau_{l_p}s\alpha]_k
    -\left[\frac{\Ad_{\e^{-sp}}-1}{\ad_{sp}}
    \tau_{l_p}(\widetilde{s\alpha}_p+\p_\l l_ps\alpha)\right]_k.
  \end{align*}
  Thus, by equation~\eqref{eq:pmdadtauzero},
  \begin{equation*}
    \p_\l\left(\frac{\Ad_{\e^{-sp}}-1}{\ad_{sp}}\tau_{l_p}s\alpha
      +\frac{\Ad_{\e^{-sp}}-1+\ad_{sp}}{\ad_{sp}^2}s\alpha
    \right)=0
  \end{equation*}
  and equation~\eqref{eq:pldadtaualpha} is proved. 
\end{proof}

The following theorem gives an explicit analytic formula for $\lcan$:
\begin{thm}\label{th:analyticity}
  The canonical $\ell$-matrix of Theorem~\ref{th:l-matrices generales}
  is analytic around $0$, and is explicitely given by:
  \begin{align}
    \label{eq:explicit l on mm}
    \p_\m l_p\xi&=(\Id_\m-R_pS_pi_\m)^{-1}R_p(S_p\xi-\xi)\\*
    \label{eq:explicit l on lm}
    \p_\l l_p\xi&=(\Id_\l-S_pR_pi_\l)^{-1}S_p(R_p\xi-\xi)\\*
    \label{eq:explicit l on ml}
    \p_\m l_ps\alpha&=(\Id_\m-R_pS_pi_\m)^{-1}R_p
    \left(S_p+K_p+\frac{\Ad_{\e^{sp}}-1-\ad_{sp}}{\ad_{sp}}\right)s\alpha\\*
    \label{eq:explicit l on ll}
    \p_\l l_ps\alpha&=-(\Id_\l-S_pR_pi_\l)^{-1}
    \left(S_p+K_p+S_pR_p\frac{\Ad_{\e^{sp}}-1-\ad_{sp}}{\ad_{sp}}\right)s\alpha
  \end{align}
  for all $\xi\,\in\l^\perp$ and $\alpha\in\l^*$, where $R$ and $S$ are
  given by:
  \begin{align}
    \label{eq:Kp}
    K_p&=\left(\p_\l\frac{\Ad_{\e^{-sp}}-1}{\ad_{sp}}i_\l\right)^{-1}\p_\l
    \frac{\Ad_{\e^{-sp}}-1+\ad_{sp}}{\ad_{sp}^2}\\*
    \label{eq:Rp}
    R_p&=(\p_\m\Ad_{\e^{-sp}}i_\m)^{-1}\p_\m\Ad_{\e^{-sp}}\\*
    \label{eq:Sp}
    S_p&=\left(\p_\l\frac{\Ad_{\e^{-sp}}-1}{\ad_{sp}}i_\l\right)^{-1}
    \p_\l\frac{\Ad_{\e^{-sp}}-1}{\ad_{sp}}.
  \end{align}
\end{thm}
\begin{proof}
  Let $\xi\in\l^\perp\subset\g^\star_p$. We get from
  equation~\eqref{eq:pmadtauzero}: $\p_\m l_p\xi=-R_p(\p_\l
  l_p\xi+\xi)$, and from equation~\eqref{eq:pmdadtauzero}: $\p_\l
  l_p\xi=-S_p(\p_\m l_p\xi+\xi)$. Combining these two equations, and the
  fact that both $(\Id_\m-R_pS_pi_\m)$ and $(\Id_\l-S_pR_pi_\l)$ are
  analytically invertible around $0$ yields equations~\eqref{eq:explicit
    l on mm} and~\eqref{eq:explicit l on lm}. 

  Let $\alpha\in\l^*$. We get from equation~\eqref{eq:pmadalpha}: $\p_\m
  l_ps\alpha=-R_p\left(\p_\l l_ps\alpha+s\alpha+
    \frac{1-\Ad_{\e^{sp}}}{\ad_{sp}}s\alpha\right)$, and from
  equation~\eqref{eq:pldadtaualpha}: $\p_\l
  l_ps\alpha=-K_ps\alpha-S_ps\alpha-S_p\p_\m l_ps\alpha$. Combining
  these two equations yields equations~\eqref{eq:explicit l on ml}
  and~\eqref{eq:explicit l on ll}.
\end{proof}

\begin{exmp}[$\ell$-matrices for a Lie quasi-bialgebra canonically
  compatible with a reductive decomposition]
  \label{ex:lmatr}We say (see~\cite{PP05}) that the Lie quasi-bialgebra
  $\G=(\g,[\,,\,],\varpi,\varphi)$ is canonically compatible with the
  reductive decomposition $\g=\l\oplus\m$ if the following conditions
  hold:
  \begin{gather}
    \varpi_\l=0\\*
    \<\m^\perp,\varpi_\bullet\m^\perp\>=0\\*
    \varphi\in\Alt(\l\otimes\l\otimes\l\,\oplus\,\l\otimes\m\otimes\m).
  \end{gather}
  In this case, the expression $S_p$ of equation~\eqref{eq:Sp} vanishes,
  and the expression $R_p$ of equation~\eqref{eq:Rp} reads as:
  \begin{equation}
    R_p=\p_\m\left(\p_\g\Ad_{\e^{-sp}}i_\g\right)^{-1}i_\g\Ad_{\e^{-sp}}.
  \end{equation}
  Thus, $\lcan$ has the simpler expression of~\cite{PP05}:
  \begin{align}
    \lcan_ps\alpha&=
    \left(\coth{\ad_{sp}-\frac1{\ad_{sp}}}\right)s\alpha,
    \ \alpha\in\l^*\\*
    \label{eq:lcanxicomp}
    \lcan_p\xi&=
    -\left(\p_\g\Ad_{\e^{-sp}}i_\g\right)^{-1}i_\g\Ad_{\e^{-sp}}\xi,
    \xi\in\l^\perp.
  \end{align}
  Morever, if $\varpi=0$, then equation~\eqref{eq:lcanxicomp} reads as:
  \begin{equation}
    \lcan_p\xi=\tanh\ad_{sp}\xi.
  \end{equation}
  This example can be specified to the Etingof--Varchenko case (see
  Example~\ref{ex:dualEV}).
\end{exmp}

\begin{exmp}[Alekseev--Meinrenken $r$-matrices~\cite{AM00}]
  Alekseev--Meinrenken $r$-matrices are obtained when $\l=\g$, for a
  cocommutative Lie quasi-bialgebra $\G=(\g,[\,,\,],0,\varphi)$. In this
  case (as $\m=\{0\}$), the only line of interest in
  Theorem~\ref{th:analyticity} is equation~\eqref{eq:explicit l on ll},
  which reads as
  \begin{equation}
    l_p\alpha=(S_p+K_p)\alpha,
  \end{equation}
  for $\alpha\in\g^*$. Also, since
  \begin{xalignat}{2}
    \ad_p\g&\subset\g^*,&\ad_p\g^*&\subset\g,
  \end{xalignat}
  the maps $K$ and $S$ read as:
  \begin{xalignat}{2}
    S_p\alpha&=\left(\frac1{\sh\ad_p}-\coth\ad_p\right)\alpha,&
    K_p\alpha&=\left(\frac1{\ad_p}-\frac1{\sh\ad_p}\right)\alpha,
  \end{xalignat}
  for all $\alpha\in\g^*$, thus we obtain:
  \begin{equation}
    \lcan_p\alpha=\rAM_p\alpha=\left(\coth\ad_p-\frac1{\ad_p}\right)\alpha,
  \end{equation}
  for all $\alpha\in\g^*$. The Alekseev-Meinrenken $r$-matrix associated
  with a cocommutative Lie quasi-bialgebra was already constructed
  in~\cite{EE03}.
\end{exmp}

\begin{exmp}[$r$-matrices for the non-compatible case]
  \label{ex:lmatrnotcomp}
  In this example, we assume that $\varpi=0$ and that
  $\varphi\in\Alt(\l\otimes\l\otimes\m)$. We set:
  \begin{xalignat}{2}
    \chpm x&=\frac12(\ch x\pm\cos x),&
    \shpm x&=\frac12(\sh x\pm\sin x),
  \end{xalignat}
  (notice that $\chpm$ and $\shpm$ are respectively the parts of order
  $0$, $2$, $1$ and $3$ modulo $4$ of the exponential map).

  One has:
  \begin{xalignat}2
    \ad_{sp}\l&\subset\m^\perp,&
    \ad_{sp}\m^\perp&\subset\m,\\*
    \ad_{sp}\m&\subset\l^\perp,&
    \ad_{sp}\l^\perp&\subset\l.
  \end{xalignat}
  Thus, for all $k\in\Nset$ such that $k\not\equiv0\mod4$ one has
  $\p_\l\ad_{sp}^ki_\l=0$ and $\p_\m\ad_{sp}^ki_\m=0$, thus the map $K$
  of Theorem~\ref{th:analyticity} reads as:
  \begin{equation}
    K_ps\alpha=\left(\frac1{\ad_{sp}}
      -\frac1{\shp\ad_{sp}}s\alpha\right),
  \end{equation}
  for all $\alpha\in\l^*$, and the maps $R$ and $S$ are represented by
  the following block-matrices (relatively to the vector space
  decomposition $\g=\l\oplus\m^\perp\oplus\m\oplus\l^\perp$, in this
  order):
  \begin{xalignat}{2}
    R_p&=
    \begin{pmatrix}
      0&0&0&0\\*
      0&0&0&0\\*
      \frac{\chm\ad_{sp}}{\chp\ad_{sp}}&
      -\frac{\shp\ad_{sp}}{\chp\ad_{sp}}&
      1&
      -\frac{\shm\ad_{sp}}{\chp\ad_{sp}}\\*
      0&0&0&0
    \end{pmatrix},&
    S_p&=
    \begin{pmatrix}
      -1&
      \frac{\chp\ad_{sp}}{\shp\ad_{sp}}&
      -\frac{\shm\ad_{sp}}{\shp\ad_{sp}}&
      \frac{\chm\ad_{sp}}{\shp\ad_{sp}}\\*
      0&0&0&0\\*
      0&0&0&0\\*
      0&0&0&0
    \end{pmatrix},
  \end{xalignat}
  and for $\alpha\in\l^*$,
  $R_p\frac{\Ad_{\e^{sp}}-1-\ad_{sp}}{\ad_{sp}}s\alpha$ reads as:
  \begin{equation}
    R_p\frac{\Ad_{\e^{sp}}-1-\ad_{sp}}{\ad_{sp}}s\alpha=
    \frac{\ad_{sp}\shp(\ad_{sp})-\chm(\ad_{sp})}
    {\ad_{sp}\chp(\ad_{sp})}s\alpha.
  \end{equation}
  Also, the mappings $(\Id_\m-R_pS_pi_\m)$ and $(\Id_\l-S_pR_pi_\l)$
  read as:
  \begin{align}
    \Id_\m-R_pS_pi_\m&=\frac{\chp\ad_{sp}\,\shp\ad_{sp}+
      \chm\ad_{sp}\,\shm\ad_{sp}}
    {\chp\ad_{sp}\,\shp\ad_{sp}}i_\m,\\*
    \Id_\l-S_pR_pi_\l&=\frac{\chp\ad_{sp}\,\shp\ad_{sp}+
      \chm\ad_{sp}\,\shm\ad_{sp}}
    {\chp\ad_{sp}\,\shp\ad_{sp}}i_\l.
  \end{align}
  We set
  \begin{align}
    F(z)&=\frac{\cosh(z)\cos(z)-1}{\cosh(z)\sin(z)+\cos(z)\sinh(z)}
    =\sum_{k\geq0}F_kz^{4k+3},\\*
    G(z)&=\frac{\sinh(z)\sin(z)}{\cosh(z)\sin(z)+\cos(z)\sinh(z)}
    =\sum_{k\geq0}G_kz^{4k+1},\\*
    H(z)&=\frac{\cos(z)\bigl(z\cosh(z)-\sinh(z)\bigr)-\sin(z)\cosh(z)+z}
    {z(\cos(z)\sinh(z)+\cosh(z)\sin(z))}
    =\sum_{k\geq0}H_kz^{4k+3}.
  \end{align}

  After computation and simplification, one obtains:
  \begin{align}
    \p_\m\lcan_p\xi&=F(\ad_{sp})\xi,\\*
    \p_\l\lcan\xi&=G(\ad_{sp})\xi,\\*
    \p_\m\lcan s\alpha&=G(\ad_{sp})s\alpha,\\*
    \p_\l\lcan s\alpha&=H(\ad_{sp})s\alpha,
  \end{align}
  for all $\xi\in\l^\perp$, $\alpha\in\l^*$.

  Equation~\eqref{eq:diffeq pour lcan} for $\lcan$ is equivalent to the
  following differential system:
  \begin{equation}
    \left\{
      \begin{aligned}
        zF'(z)&=-z\bigl(F(z)^2+G(z)^2\bigr)\\*
        zG'(z)&=z-zG(z)\bigl(H(z)+F(z)\bigr)-G(z)\\*
        zH'(z)&=-zG(z)^2-zH(z)^2-2H(z).
      \end{aligned}
    \right.
  \end{equation}
\end{exmp}

\subsection{Functoriality}
\begin{defn}
  A \emph{bidynamical Lie quasi-bialgebra (over $\l$)} is a Lie
  quasi-bialgebra $\G=(\g,[\,,\,],\varpi,\varphi)$ such that
  $\varpi_\l=0$, $\varphi\equiv0\mod\l$ and such that there exists a
  reductive decomposition $\g=\l\oplus\m$.

  A morphism $\phi$ between two bidynamical Lie quasi-bialgebras
  over $\l$, $\G=(\g=\l\oplus\m,[\,,\,],\varpi,\varphi)$ and
  $\G'=(\g'=\l\oplus\m',[\,,\,]',\varpi',\varphi')$, is a Lie algebra
  morphism $\phi\colon\g\to\g'$ such that $\phi(z)=z$ for all $z\in\l$,
  $\phi(\m)\subset\m'$ and
  \begin{align}
    \phi\varpi_x\phi^*&=\varpi'_{\phi(x)},\quad\forall x\in\g,\\*
    \phi^{(3)}\varphi&=\varphi'.
  \end{align}
\end{defn}

\begin{prop}
  \label{pr:upsilvfnqerjkl}
  Let $\G=(\g=\l\oplus\m,[\,,\,],\varpi,\varphi)$ and
  $\G'=(\g'=\l\oplus\m,[\,,\,]',\varpi',\varphi')$ be two bidynamical
  Lie quasi-bialgebras, and let $\upsi\colon\G\to\G'$ be a bidynamical
  Lie quasi-bialgebra morphism. Then the following equality holds
  \begin{equation}
    \upsi\lcanG{\G,\l,\m}\upsi^*=\lcanG{\G',\l,\m'}.
  \end{equation}
  In particular, the Lie groupoid morphism
  \begin{equation*}
    \appli\Psi{\GG=U\times G\times U}{\GG'=U\times G'\times U}
    {(p,x,q)}{(p,\psi(x),q)}
  \end{equation*}
  is a Poisson groupoid morphism, when $\GG$ and $\GG'$ are equipped
  with the Poisson bracket induced by $\lcanG{\G,\l,\m}$ and
  $\lcanG{\G',\l,\m'}$ respectively, where $\psi\colon G\to G'$ is the
  Lie group morphism integrating the Lie algebra morphism $\upsi$. 
\end{prop}
\begin{proof}
  We set $\lambda_p=\upsi\lcanG{\G,\l,\m}\upsi^*$. Then, it is easy to
  check, from equation~\eqref{eq:diffeq pour lcan}, that $\lambda$
  satisfies the same equation as $\lcanG{\G',\l,\m'}$, namely:
  \begin{equation*}
    \d_p\lambda(p)=\tau_{-\lambda_p}\ad'_{s'p}\tau_{\lambda_p}
    -\p'_\l\lambda_p-\lambda_ps'i^{\prime*},
  \end{equation*}
  thus $\upsi\lcanG{\G,\l,\m}\upsi^*=\lcanG{\G',\l,\m'}$.
\end{proof}

\section{Trivialization and duality}
\label{sec:TandD}
From now on, we denote by $U$ the domain of analyticity of the canonical
$\ell$-matrix $\lcan$. Notice that $U$ is $\l$-equivariant, but not
simply-connected in general.

\subsection{Trivial Lie algebroids}
Let $(\g,[\,,\,]_\g)$ be a Lie algebra and $M$ a manifold. Recall
(see~\cite{M87}) that the trivial Lie algebroid on $M$ with vertex
algebra $\g$ is the vector bundle $A=\T M\oplus(M\times\g)=\T M\times\g$
over $M$ (Whitney sum), where the anchor is the projection on $\T M$,
and the bracket is defined as follows: let $\sigma$ and $\sigma'$ be two
sections of the vector bundle $A$, say $\sigma=(X,x)$ and
$\sigma'=(X',x')$ where $X$ and $X'$ are two vector fields on $M$ and
$x,\,x'\colon M\to\g$, and set
\begin{equation}\label{eq:trivial algebds bracket}
  [\sigma,\sigma']_A=[X,X']\oplus
  \left(X\cdot x'-X'\cdot x+[x,x']_\g\right)
\end{equation}
The bracket in the first component of the right hand side of
equation~\eqref{eq:trivial algebds bracket} is the bracket of vector
fields on $M$, and $X\cdot x'$ denotes the derivative of $x'$ in the
direction of $X$.

The Lie algebroid of the trivial groupoid $\GG=U\times G\times U$ is the
trivial Lie algebroid over $U$: $\Alg(\GG)=U\times(\l^*\oplus\g)$, with
the following bracket on its sections:
\begin{equation}
  [\sigma,\sigma']_p=\bigl(\d_p\alpha'(\alpha_p)-\d_p\alpha(\alpha'_p),
  \d_px'(\alpha_p)-\d_px(\alpha'_p)+[x_p,x'_p]_\g\bigr),
\end{equation}
for $p\in U$ and $\sigma_p=\alpha_p+x_p$,
$\sigma'_p=\alpha'_p+x'_p\in\l^*\oplus\g$, and the anchor is:
\begin{equation}
  a^{\Alg(\GG)}\sigma=\alpha.
\end{equation}

\subsection{Trivialization}
We recall (see~\cite{LP05}) that the Lie algebroid of the dual of the
Poisson groupoid $\GG$ is the vector bundle $N(U)=U\times\l\times\g^*$
over $U$, together with the following bracket on its sections:
\begin{equation}
\begin{aligned}
  \big[(z,\xi),(z',\xi')\big]^{N(U)}_p=
  &\Big(\d_pz'\big(a^{N(U)}_p(z_p,\xi_p)\big)
    -\d_pz\big(a^{N(U)}_p(z'_p,\xi'_p)\big)-[z_p,z'_p]+
    \<\xi,\d_pl(\cdot)\xi'\>,\\*
  &\qquad
    \d_p\xi'\big(a^{N(U)}_p(z_p,\xi_p)\big)
    -\d_p\xi\big(a^{N(U)}_p(z'_p,\xi'_p)\big)
  +\ad^*_{z_p}\xi'_p-\ad^*_{z_p'}\xi_p\\*
  &\qquad+\<\xi_p,\varpi_\bullet\xi'_p\>
  +\ad^*_{l_p\xi_p}\xi'_p-\ad^*_{l_p\xi'_p}\xi_p\Big)
\end{aligned}
\end{equation}
and the anchor:
\begin{equation}
  a^{N(U)}_p(z,\xi)=i^*\xi-\ad^*_zp.
\end{equation}

We want a trivialization of the Lie algebroid $N(U)$, that is a Lie
algebroid isomorphism $T\colon\Alg(\GG^\star)\to N(U)$, where
$\GG^\star=U\times G^\star_0\times U$ is the trivial groupoid over $U$,
with vertex group $G^\star_0$ such that $\Lie(G^\star_0)=\g^\star_0$.
Such an isomorphism can be split into two parts (see~\cite{M87}), a Lie
algebra bundle isomorphism $\psi\colon U\times\g_0^\star\to\Ker
a^{N(U)}$ and a flat connection $\nabla\colon U\times\l^*\to N(U)$ such
that $[\nabla(\alpha),\psi X]^{N(U)}=\psi(\d X(\alpha))$ for any smooth
section $\alpha\in\Gamma(U\times\l^*)$ and
$X\in\Gamma(U\times\g^\star_0)$ --- we recall that by definition,
$\nabla$ is a flat connection if $[\nabla\alpha,\nabla\beta]=0$ for all
$\alpha,\,\beta\in\Gamma(U\times\l)$. Then, setting
$T(\alpha+X)=\nabla\alpha+\psi X$ for $\alpha\in\l^*$, $X\in\g^\star_0$
provides a trivialization $T$.

We denote by $U$ the domain of analycity of $\lcan$ (which is
$\Ad^*_L$-equivariant since $\lcan$ is $\l$-equivariant and $L$ is
connected). As a corollary of Theorem~\ref{th:analyticity} we have the
following:
\begin{cor}\label{cor:LAB morph}
  Let $l=\lcanG{\G,\l,\m}\in\Dynl_0(U,\G)$. Then, for all
  $X_p\in\g^\star_p$, the expression
  \begin{equation}
    \label{eq:expressioningstar0}
    \Ad_{\e^{-sp}}\tau_{l_p}X_p
  \end{equation}
  lies in $\g^\star_0$. In particular, the map
  \begin{equation*}
    \appli{\phi_p}{\g^\star_p}{\g^\star_0}{X_p}
    {\Ad_{\e^{-sp}}\tau_{l_p}X_p}
  \end{equation*}
  is a Lie algebra isomorphism, and its inverse is given by
  \begin{equation*}
    \appli{\phi^{-1}_p}{\g^\star_0}{\g^\star_p}{z+\xi}
    {\left(\p_\l\dfrac{\Ad_{\e^{sp}-1}}{\ad_{sp}}
        +\p_{\g^*}\Ad_{\e^{sp}}\right)(z+\xi).}
  \end{equation*}
\end{cor}
\begin{proof}
  By equation~\eqref{eq:pmadtauzero}, we know that the
  expression~\eqref{eq:expressioningstar0} lies in
  $\l\oplus\g^*\subset\dlie$. Also, since
  $\p_{\m^\perp}\ad_{sp}=\ad_{sp}\p_\l$, applying $\ad_{sp}$ to both
  sides of equation~\eqref{eq:pmdadtauzero} yields:
  \begin{equation*}
    \p_{\m^\perp}\Ad_{\e^{-sp}}\tau_{l_p}X_p-\p_{\m^\perp}X_p=-\ad_{sp}\p_\l
    X_p.
  \end{equation*}
  Thus, since $\p_{\m^\perp}X_p=-\ad_{sp}\p_\l$, the expression
  $\p_{\m^\perp}\Ad_{\e^{-sp}}\tau_{l_p}X_p$ vanishes, and
  expression~\eqref{eq:expressioningstar0} lies in
  $\l\oplus\l^\perp=\g_0^\star$. Clearly, the map $\phi_p$ is a Lie
  algebra isomorphism. A simple computation using the relations of
  Proposition~\ref{pr:pmadtauzero} shows that
  $\left(\p_\l\frac{\Ad_{\e^{sp}-1}}{\ad_{sp}}
    +\p_{\g^*}\Ad_{\e^{sp}}\right)\Ad_{\e^{-sp}}
  \tau_{l_p}(z+\ad_{sp}z+\xi)=z+\ad_{sp}z+\xi$ for all $z\in\l$ and
  $\xi\in\l^\perp$.
\end{proof}
Corollary~\ref{cor:LAB morph} implies that the bundle map
\begin{equation}
  \appli{\psi}{U\times\g^\star_0}
  {\Ker a^{N(U)}\subset U\times(\l\oplus\g^*)}
  {(p,X)}{\left(p,-\phi_p^{-1}X\right)}
\end{equation}
is a Lie algebra bundle isomorphism. So, in order to complete the
trivialization, we need a flat connection $\nabla\colon U\times\l^*\to
N(U)$, satisfying $[\nabla(\alpha),\psi X]^{N(U)}=\psi(\d X(\alpha))$
for any smooth section $\alpha\in\Gamma(U\times\l^*)$ and
$X\in\Gamma(U\times\g^\star_0)$. 

For all $\alpha\in\l^*$, we set
\begin{equation}\label{eq:form de nabla}
  \nabla_p(\alpha)=(u_p\alpha,s\alpha+\ad_{sp}u_p\alpha+v_p\alpha),
\end{equation}
where $u_p\colon\l^*\to\l$ and $v_p\colon\l^*\to\l^\perp$ are given by:
\begin{equation}\label{eq:sol uv}
  \left\{
    \begin{aligned}
      u_p\alpha&=\p_\l\tau_{-l_p}\frac{\Ad_{\e^{sp}}-1}{\ad_{sp}}s\alpha
      =\p_\l\frac{\Ad_{\e^{sp}}-1-\ad_{sp}}{\ad_{sp}^2}s\alpha\\*
      v_p\alpha&=\p_{\l^\perp}\frac{\Ad_{\e^{sp}}-1}{\ad_{sp}}s\alpha
    \end{aligned}
    \right. 
\end{equation}
for all $\alpha\in\l^*$. We show in Proposition~\ref{pr:nablaflat} below
that $\nabla$ is a flat connection.

We introduce the following notations:
\begin{align}
  \happy\alpha&=\p_\g\frac{\Ad_{\e^{sp}}-1}{\ad_{sp}}s\alpha\\*
  \saddy\alpha&=\p_{\g^*}\frac{\Ad_{\e^{sp}}-1}{\ad_{sp}}s\alpha
\end{align}
for $\alpha\in\l^*$. Then, $\nabla$ is given by:
\begin{equation}
  \nabla_p\alpha=\left(u_p\alpha,\saddy\alpha\right)
\end{equation}
for all $\alpha\in\l^*$. The following lemma, the proof of which is
straightforward, will help in our computations:
\begin{lem}\label{lm:u=h+s}
  For all $\alpha\in\l^*$,
  \begin{align}
    u_p\alpha&=\happy\alpha-l_p\saddy\alpha\\*
    &=\p_\g\tau_{-l_p}\frac{\Ad_{\e^{sp}}-1}{\ad_{sp}}s\alpha.
  \end{align}
\end{lem}

\begin{prop}\label{pr:nablaflat}
  The map $\nabla$ is a flat connection.
\end{prop}

\begin{proof}
  To show that $\nabla$ is a flat connection, we have to show the two
  following equalities for all $\alpha,\,\beta\in\l^*$ seen as constant
  sections:
  \begin{equation}
    \label{eq:flat1}
    \d_pu(\alpha)\beta-\d_pu(\beta)\alpha-[u_p\alpha,u_p\beta]
    +\<\saddy\alpha,\d_pl(\cdot)\saddy\beta\>=0
  \end{equation}
  \begin{multline}
  \label{eq:flat2}
    \d_p\saddyk\beta(\alpha)-\d_p\saddyk\alpha(\beta)
    -\p_{\g^*}[u_p\alpha,\saddy\beta]-\p_{\g^*}[\saddy\alpha,u_p\beta]
    -\p_{\g^*}[\saddy\alpha,\saddy\beta]\\*
    -\p_{\g^*}[l_p\saddy\alpha,\saddy\beta]
    -\p_{\g^*}[\saddy\alpha,l_p\saddy\beta]=0.
  \end{multline}
  
  The left hand side of~\eqref{eq:flat2} vanishes, since
  $\d_p\saddyk\beta(\alpha)-\d_p\saddyk\alpha(\beta)=
  \left[\happy\alpha+\saddy\alpha,\happy\beta+\saddy\beta\right]$, and
  by Lemma~\ref{lm:u=h+s}.

  Equation~\eqref{eq:flat1} vanishes too, by~\eqref{eq:CDYBEg}
  and~\eqref{eq:l-eqg}.
\end{proof}

The trivialization is given in the following theorem:
\begin{thm}
  \label{th:bundleisom}
  Let $l$ be a canonical dynamical $\ell$-matrix. Then, the bundle map:
  \begin{equation}
    T\colon U\times\l^*\times\g^\star_0\longrightarrow
    N(U)=U\times\l\times\g^*
  \end{equation}
  given by
  \begin{equation}\label{eq:trivialisation}
    T_p(\alpha,z+\xi)=
    \left(u_p\alpha-\p_\l\tau_{-l_p}\Ad_{\e^{sp}}(z+\xi),
      \saddy\alpha-\p_{\g^*}\Ad_{\e^{sp}}(z+\xi)\right)
  \end{equation}
  is a Lie algebroid isomorphism. 
\end{thm}
\begin{proof}
  It only remains to show that $[\nabla(\alpha),\psi X]^{N(U)}=0$, for
  all $\alpha\in\l^*$, and all $X\in\g^\star_0$ seen as constant
  sections. 
  
  First, notice that $\psi_pz=-\nabla_p\ad^*_zp-(z,0)$ for all
  $z\in\l$. Thus, $[\nabla\alpha,\psi z]_p^{N(U)}=
  -\nabla_p[\alpha,\ad^*_z\cdot]_p-[\nabla\alpha,(z,0)]^{N(U)}_p=
  \psi_p\bigl(\d_pz(\alpha)\bigr)$.
  
  Now, if $\xi\in\l^\perp$ (seen as a constant section), the first
  component of
  $[\nabla\alpha,\psi\xi]^{N(U)}_p$ is:
  \begin{align*}
    \p_\l\d_pl(\alpha)\Ad_{\e^{sp}}\xi
    &-\p_\l\tau_{-l_p}\left[\frac{\Ad_{\e^{sp}}-1}{\ad_{sp}}s\alpha,
      \Ad_{\e^{sp}}\xi\right]+[u_p\alpha,
    \p_\l\tau_{-l_p}\Ad_{\e^{sp}}\xi]-\<\saddy\alpha,\d_pl(\cdot)
    \p_\g^*\Ad_{\e^{sp}}\xi\>\\*
    &=
    -\p_\l\d_pl(\ad^*_{u_p\alpha}p)\p_{\g^*}\Ad_{\e^{sp}}\xi
    +\p_\l\d_pl(\ad^*_{\p_\l\tau_{-l_p}\Ad_{\e^{sp}}\xi}p)\saddy\alpha\\*
    &\qquad-\p_\l\left[\tau_{-l_p}\frac{\Ad_{\e^{sp}}-1}{\ad_{sp}}s\alpha,
      \tau_{-l_p}\Ad_{\e^{sp}}\xi\right]_{\dlie^{l_p}}
    +[u_p\alpha,\p_\g\tau_{-l_p}\Ad_{\e^{sp}}\xi]
    \\*&=0
  \end{align*}
and the second component vanishes too. 
\end{proof}

The trivialization may also be written:
\begin{equation}
  T_p(\alpha,z+\xi)=\tau_{-l_p}\frac{\Ad_{\e^{sp}}-1}{\ad_{sp}}s\alpha
  -\tau_{-l_p}\Ad_{\e^{sp}}(z+\xi)
\end{equation}
where $\alpha\in\l^*$, $z\in\l$ and $\xi\in\l^\perp$, and also in a way
where $l$ does not appear:
\begin{equation}
  \label{eq:Trivialwithoutl}
  T_p(\alpha,z+\xi)=\p_\l\frac{\Ad_{\e^{sp}}-1-\ad_{sp}}{\ad_{sp}^2}s\alpha
  +\p_{\g^*}\frac{\Ad_{\e^{sp}}-1}{\ad_{sp}}s\alpha
  -\p_\l\frac{\Ad_{\e^{sp}}-1}{\ad_{sp}}(z+\xi)
  -\p_{\g^*}\Ad_{\e^{sp}}(z+\xi).
\end{equation}
To show this last equality, compute the adjoint $T^*_p$ of $T_p$, and
use Proposition~\ref{pr:pmadtauzero}. 

\subsection{Duality}
For a Manin quasi-triple $(\dlie,\g,\h)$ we denote by
$\G_{(\dlie,\g,\h)}$ the corresponding Lie quasi-bialgebra; for a Lie
quasi-bialgebra $\G=(\g,[\,,\,],\varpi,\varphi)$, we set
$\G^-=(\g,[\,,\,],-\varpi,\varphi)$.

We start with the definition of a duality for bidynamical Lie
quasi-bialgebras, which was introduced in~\cite{PP05}.

\begin{defn}\label{df:q-big dual}
  Let $\G=(\g=\l\oplus\m,[\,,\,],\varpi,\varphi)$ be a bidynamical Lie
  quasi-bialgebra with canonical double $\dlie$. The Lie quasi-bialgebra
  \begin{equation*}
    \G^\star=\left(\G_{(\dlie,\l\oplus\l^\perp,\m^\perp\oplus\m)}\right)^-
  \end{equation*}
  is called the \emph{dual over $\l$ of the bidynamical Lie
    quasi-bialgebra $\G$}.
\end{defn}
Observe that $\g^\star=\l\oplus\l^\perp$ is indeed a lagrangian
subalgebra of $\dlie$, so that the dual over $\l$ is well-defined. Also
observe that the dual of a bidynamical Lie quasi-bialgebra is again a
bidynamical Lie quasi-bialgebra.

Let $\op\colon\g\to\g$ be the standard involution associated with the
reductive decomposition $\g=\l\oplus\m$:
\begin{equation}
  \op(z)=z\qquad\qquad\op(u)=-u
\end{equation}
for all $z\in\l$ and $u\in\m$. We define the Lie bracket $[\,,\,]^\op$
on the vector space $\g$ as follows:
\begin{xalignat}{3}
  [z,z']^\op&=[z,z']&
  [z,u]^\op&=-[z,u]&
  [u,u']^\op&=[u,u']
\end{xalignat}
for all $z,\,z'\in\l$, $u,\,u'\in\m$, and we denote by $\g^\op$ the
resulting Lie algebra. We also set
$\G^\op=(\g^\op,[\,,\,]^\op,\varpi^\op,\varphi^\op)$ where
\begin{xalignat}{2}
  \varpi^\op_x&=\op\,\varpi_{\op(x)}\,\op^*&
  \varphi^\op&=(\op\otimes\op\otimes\op)\,\varphi
\end{xalignat}

We denote by $\dlie^\star$ the canonical double of $\G^\star$. First,
observe that under the canonical vector space identification
$\dlie^\star\simeq\dlie$, the Lie algebra $\g$ is not a Lie subalgebra
of $\dlie^\star$ (but the Lie algebra $\g$ is isomorphic to the Lie
subalgebra $\g^\op=\l\oplus\m$ of $\dlie^\star$).  Second, observe that
under the canonical identification $\dlie^\star\simeq\dlie$, then
$(\G^\star)^\star\not=\G$, but rather $(\G^\star)^\star=\G^\op$, which
is isomorphic to $\G$.

We now turn to our main duality statement which provides the dual
Poisson groupoid of a Poisson groupoid associated with a canonical
$\ell$-matrix:
\begin{thm}\label{th:duality}
  The dual Poisson groupoid of the dynamical Poisson groupoid associated
  with a canonical $\ell$-matrix $\lcanG{\G,\l,\m}$ is isomorphic to (a
  covering of) the dynamical Poisson groupoid $U\times G^\star\times U$
  with the Poisson structure associated with the canonical $\ell$-matrix
  on $U$ $\lcanG{\G^\star,\l,\l^\perp}$, where $G^\star$ is the
  connected, simply connected Lie group with Lie algebra $\g^\star$. 
\end{thm}
\begin{proof}
  We compute $-T^*\colon U\times\l^*\times\g\to
  U\times\l\times(\g^\star_0)^*$ from
  equation~\eqref{eq:Trivialwithoutl}.
  \begin{multline}
    -T_p^*(\alpha,x)=\twolinesleft(
    -\p_\l\frac{\Ad_{\e^{-sp}}-1+\ad_{sp}}{\ad_{sp}^2}s\alpha
    +\p_\l\frac{\Ad_{\e^{sp}}-1}{\ad_{sp}}x,\twolinesmid
    -(\p_{\m^\perp}+\p_{\m})\frac{\Ad_{\e^{-sp}}-1}{\ad_{sp}}s\alpha
    +(\p_{\m^\perp}+\p_{\m})\Ad_{\e^{-sp}}x\twolinesright)
  \end{multline}
  for $\alpha\in\l^*$ and $x\in\g$. 

  Now let $T^\star$ be the algebroid isomorphism
  \begin{equation*}
    T^\star\colon U\times\l^*\times\g^\op\longrightarrow
    U\times\l\times\bigl(\g^\star_0\bigr)^*
  \end{equation*}
  given by Theorem~\ref{th:bundleisom}, associated with the datum
  $\G^\star$. An easy computation shows that
  $T^\star\circ\widehat\op=-T^*$, \ie that $-T^*$ is indeed a Lie
  algebroid isomorphism, where $\widehat\op\colon U\times\l^*\times\g\to
  U\times\l^*\times\g^\op$ is the Lie algebroid isomorphism given by
  $\widehat\op_p(\alpha,x)=\bigl(\alpha,\op(x)\bigr)$ for all $p\in U$,
  $\alpha\in\l^*$ and $x\in\g$.
\end{proof}

In this part, we also showed the following:
\begin{thm}
  A dynamical Poisson groupoid on $U\subset\l^*$ with $0\in U$ is
  bidynamical if and only if there is a reductive decomposition
  $\g=\l\oplus\m$ of the Lie algebra $\g$ of the vertex group.
\end{thm}
\begin{proof}
  It is shown in~\cite{PP05} that it is a necessary condition.
  Theorem~\ref{th:duality} shows that it is indeed a sufficient
  condition.
\end{proof}

\begin{exmp}[Dual of the Alekseev-Meinrenken $r$-matrix]
  Let $\G=(\g,[\,,\,],0,\varphi)$ be a cocommutative Lie
  quasi-bialgebra. Clearly, the dual over $\g$ of $\G$ is $\G$ itself,
  so that $\rAM$ is self-dual.
\end{exmp}

\begin{exmp}[Dual $\ell$-matrix of the Etingof--Varchenko $r$-matrices]
  \label{ex:dualEV}
  Let $\g$ be a semi-simple Lie algebra, $\l$ a Cartan subalgebra of
  $\g$, $\Delta$ the set of roots, and $\G=(\g,[\,,\,],0,\varphi)$ the
  Lie cocommutative quasi-bialgebra with $\varphi$ given by
  $\<\xi\otimes\eta\otimes\zeta,\varphi\>=\<\xi,[\eta,\zeta]\>$ (where
  we have identified $\g\simeq\g^*$ via the Killing form). We recall
  (see~\cite{EV98}) the form of $r$-matrices in this case: let
  $\Delta^s\subset\Delta$ be a choice of simple roots. We denote by
  $\Delta^\pm$ the set of positive/negative roots in $\Delta$. Let
  $\Gamma\subset\Delta^s$, and denote by $\<\Gamma\>$ the set of roots
  spanned by $\Gamma$, and set
  $\bar\Gamma^\pm=\Delta^\pm\setminus\<\Gamma\>$. Let $\mu\in\h^*$. For
  all root $\alpha\in\Delta$ and $p\in\h^*$ we define
  \begin{equation*}
    t_\alpha(p)=
    \begin{cases}
      \coth\bigl((\alpha,p+\mu)_{\h^*}\bigr)&\text{if
        $\alpha\in\<\Gamma\>$,}\\*
      \pm1&\text{if $\alpha\in\bar\Gamma^\pm$}
    \end{cases}
  \end{equation*}
  ($\mu$ is assumed to be chosen such that $t_\alpha$ is defined for
  $p=0$). The Etingof--Varchenko dynamical $r$-matrix associated with
  $\Gamma$ is given as follows:
  \begin{equation}
    \rEV_pe^\alpha=t_\alpha(p)e_{-\alpha},\quad
    \rEV_p\alpha=0,\ \alpha\in\Delta.
  \end{equation}
  We set $\t=\rEV_0\in\bigwedge^2\g$, and $t_\alpha=t_\alpha(0)$. Using
  Proposition~\ref{pr:Dynltwist}, we see that $\rEV-\t$ is a dynamical
  $\ell$-matrix for the Lie quasi-bialgebra
  $\G^\t=(\g,[\,,\,],\t,\varphi^\t)$. The dual over $\l$ of $\G^\t$ is
  $\G^\star=(\g^\star,[\,,\,]_{\g^\star},\varpi^\star,\varphi^\star)$
  with
  \begin{xalignat*}{2}
    \varpi^\star_{e^\alpha}\beta&=-t_\alpha(\alpha,\beta)e^\alpha,&
    \<\alpha\otimes\beta\otimes1,\varphi^\star\>&=0,\\*
    \varpi^\star_{e^\alpha}e_\gamma&=
    -N_{\gamma,\alpha-\gamma}e^{\alpha-\gamma},&
    \<\alpha\otimes e_\beta\otimes1,\varphi^\star\>&=
    (\alpha,\beta)e^{-\beta},\\*
    \varpi^\star_{e^\alpha}e_\alpha&=t_\alpha h_\alpha,&
    \<\e_\alpha\otimes
    e_\beta\otimes1,\varphi^\star\>&=\delta_{\alpha,-\beta}h_\alpha,
  \end{xalignat*}
  for $\alpha,\,\beta,\,\gamma\in\Delta$, $\alpha\neq\gamma$. We compute
  $\ad^\star_{sp}$ and its powers:
  \begin{xalignat*}2
    {\ad^\star_{sp}}^{2k}e_\alpha&=(p,\alpha)^{2k}e_\alpha,&
    {\ad^\star_{sp}}^{2k}e^\alpha&=(p,\alpha)^{2k}e^\alpha,\\*[-1.5em]
  \end{xalignat*}
  \begin{equation*}
    {\ad^\star_{sp}}^{2k+1}e_\alpha=(p,\alpha)^{2k+1}(-t_\alpha
    e_\alpha+e^{-\alpha}),
  \end{equation*}
  \begin{equation*}
    {\ad^\star_{sp}}^{2k+1}e^\alpha=
    (p,\alpha)^{2k+1}\bigl((t_\alpha^2-1)e^\alpha+t_\alpha
    e^\alpha\bigr),
  \end{equation*}
  and we obtain the form of the canonical $\ell$-matrix associated with
  $\G^\star$:
  \begin{align}
    \lEV_pe_\alpha&=
    \frac1{\coth\bigl((\alpha,p)_{\h^*}\bigr)-t_\alpha}e^{-\alpha},\\*
    \lEV_p\alpha&=0.
  \end{align}
  We end this example with the following remark: the $r$-matrix $\rEV$
  is linked to the function $x\mapsto\coth(x-a)$ which satisfies the
  following differential equation:
  \begin{equation*}
    f'(x)-f(x)^2=1,
  \end{equation*}
  and the $\ell$-matrix $\lEV$ is related to the function
  $x\mapsto\dfrac1{\coth(x)-a}$ which satisfies the following
  differential equation:
  \begin{equation*}
    f'(x)+(a^2-1)f(x)^2+2af(x)=-1.
  \end{equation*}
\end{exmp}

\begin{exmp}[Dual of the non-compatible $r$-matrix of
  Example~\ref{ex:lmatrnotcomp}]
  Let $\G=(\g=\l\oplus\m,[\,,\,],0,\varphi)$ be the bidynamical Lie
  quasi-bialgebra of Example~\ref{ex:lmatrnotcomp} (with
  $\varphi\in\Alt(\l\otimes\l\otimes\m)$). The dual over $\l$ of the
  $\G$ is $\G^\star=(\g^\star,[\,,\,]^\star,\varpi^\star,\varphi^\star)$
  where the Lie bracket $[\,,\,]^\star$ on $\g^\star$ is given by:
  \begin{xalignat}{3}
    [z,z']^\star&=[z,z']_\g,&
    [z,\xi]^\star&=-\ad_z^*\xi,&
    [\xi,\xi']^\star&=0,
  \end{xalignat}
  for all $z,\,z'\in\l$, $\xi,\,\xi'\in\l^\perp$, the cobracket
  $\varpi^\star$ is given by:
  \begin{equation}
    \varpi^\star_{z+\xi}(s\alpha+u)=
    -\<\xi\otimes s\alpha\otimes1,\varphi\>
    -\p_{\l^\perp}\ad^*_u\xi,
  \end{equation}
  for all $z\in\l$, $\alpha\in\l^*$, $u\in\m$ and $\xi\in\l^\perp$, and
  $\varphi^\star$ is given by:
  \begin{equation}
    \<(s\alpha+u)\otimes(s\beta+v)\otimes1,\varphi^\star\>=
    \p_\l[u,v]_\g+\ad^*_vs\alpha-\ad^*_us\beta,
  \end{equation}
  for all $\alpha,\,\beta\in\l^*$ and $u,\,v\in\m$. We set
  $\klie=\l\subset\g^\star$ and $\ulie=\l^\perp\subset\g^\star$, so that
  $\g^\star=\klie\oplus\ulie$. We denote by $\ads$ the adjoint action of
  the double $\dlie^\star$ of the Lie quasi-bialgebra $\G^\star$. One
  has:
  \begin{xalignat}{2}
    \ads_{sp}\klie&\subset\ulie^\perp\subset{\g^\star}^*,&
    \ads_{sp}\ulie^\perp&\subset\klie^\perp\subset{\g^\star}^*,\\*
    \ads_{sp}\klie^\perp&\subset\ulie\subset\g^\star,&
    \ads_{sp}\ulie&\subset\klie\subset\g^\star.
  \end{xalignat}
  We define the following functions:
  \begin{align}
    F^\star(z)&=\frac{2\sin(z)\sinh(z)}{\cosh(z)\sin(z)+\cos(z)\sinh(z)}
    =\sum_{k\geq0}F^\star_nz^{4k+1},\\*
    G^\star(z)&=-2\frac{\bigl(-\cos(z)\sinh(z)+\sin(z)\cosh(z)\bigr)^2}
    {\cosh(2z)\cos(2z)-1}
    =\sum_{k\geq0}G^\star_nz^{4k+2},\\*
    H^\star(z)&=\frac{2z\cos(z)\cosh(z)-\sin(z)\cosh(z)-\cos(z)\sinh(z)}
    {z(\cos(z)\sinh(z)+\sin(z)\cosh(z))}
    =\sum_{k\geq0}H^\star_nz^{4k+3}.
  \end{align}
  A computation shows that $\lcanG{\G^\star}$ reads as:
  \begin{align}
    \p_\ulie\lcanG{\G^\star}_pu&=F^\star(\ads_{sp})u,\\*
    \p_\klie\lcanG{\G^\star}_pu&=G^\star(\ads_{sp})u,\\*
    \p_\ulie\lcanG{\G^\star}_ps\alpha&=-G^\star(\ads_{sp})s\alpha,\\*
    \p_\klie\lcanG{\G^\star}_ps\alpha&=H^\star(\ads_{sp})s\alpha,
  \end{align}
  where $u\in\klie^\perp$, $\alpha\in\klie^\star$.
  Equation~\eqref{eq:diffeq pour lcan} for $\lcanG{\G^\star}$ is
  equivalent to the following the differential system:
  \begin{equation}
    \left\{
      \begin{aligned}
        z{F^\star}'(z)&=z\bigl(1+G^\star(z)^2\bigr),\\*
        z{G^\star}'(z)&=zF^\star(z)-zG^\star(z)H^\star(z)-G^\star(z),\\*
        z{H^\star}'(z)&=-2zG^\star(z)-zH^\star(z)^2-2H^\star(z).
      \end{aligned}
    \right.
  \end{equation}
\end{exmp}

\section{Link between $\lcanG{\G,\l,\m}$ and
  $\lcanG{\G^\star,\l,\l^\perp}$}
\label{sc:Link duality}

Let $\G=(\g,[\,,\,],\varpi,\varphi)$ be a Lie quasi-bialgebra. It is
well-known (see e.g.~\cite{AKS00}) that the canonical double $\dlie$ of
$\G$ carries a Lie quasi-bialgebra structure:
\begin{equation}
  \G^{(2)}=(\dlie,[\,,\,]_\dlie,\partial\rmx,j^{(3)}\varphi)
\end{equation}
where $j\colon\g\to\dlie$ is the canonical inclusion, and $\partial\rmx$
is the coboundary of the ``$r$-matrix''
$\rmx=\frac12(\p_{\g^*}-\p_\g)\in\bigwedge^2\dlie$:
\begin{equation}
  \partial_X\rmx=\ad_X\rmx+\rmx\ad^*_X,\quad X\in\dlie.
\end{equation}
It is also well-known that twisting $\G^{(2)}$ via $-\rmx$ yields the
cocommutative Lie quasi-bialgebra $\bigl(\G^{(2)}\bigr)^{-\rmx}=
(\dlie,[\,,\,],0,\frac14\<\Omega^\dlie,\Omega^\dlie\>)$ where
$\Omega^\dlie\in\Sym^2\dlie$ is the symmetric $2$-tensor associated with
the canonical bilinear form on $\dlie$.

Now, if $\G$ is a canonical dynamical Lie quasi-bialgebra, then so is
$\G^{(2)}$ with reductive decomposition $\dlie=\l\oplus(\m\oplus\g^*)$
over $\l$, and $j$ is thus a canonical dynamical Lie quasi-bialgebra
morphism.  Thus, it follows from Proposition~\ref{pr:upsilvfnqerjkl}
that
\begin{equation}
  j\lcanG{\G,\l,\m}j^*=\lcanG{\G^{(2)},\l,\m\oplus\g^*}.
\end{equation}

We denote by $K$ the Lie algebra isomorphism from the double $\dlie$ of
$\G$ to the double $\dlie^\star$ of the dual $\G^\star$ given by
$K=\p_{\l\oplus\l^\perp}-\p_{\m^\perp\oplus\m}$ when the vector spaces
$\dlie$ and $\dlie^\star$ are canonically identified. Let
$\rmx^\star=\frac12(\p_{{\g^\star}^*}-\p_{\g^\star})$ be the $r$-matrix
associated to the Manin triple $(\dlie^\star,\g^\star,\g^{\star*})$.
Clearly, $\Omega^{\dlie^\star}=(K\otimes K)\Omega^\dlie$, and
$K\colon\bigl(\G^{(2)}\bigr)^{-\rmx}\to\bigl(\G^{\star(2)}\bigr)^{-\rmx^\star}$
is thus a canonical dynamical Lie quasi-algebra morphism. Hence, using
Proposition~\ref{pr:Dynltwist} we have:
\begin{equation}
  K\left(j^\star\lcanG{\G^\star,\l,\l^\perp}
    {j^\star}^*+\rmx^\star\right)K^*
  \in\Dynl\left(U,\bigl(\G^{(2)}\bigr)^{-\rmx}\right)
\end{equation}
and thus
\begin{equation}
  \label{eq:cqzejkbvfhqzeuiobvf}
  Kj^\star\lcanG{\G^\star,\l,\l^\perp}{j^\star}^*K^*+K\rmx^\star K^*+\rmx
  \in\Dynl\bigl(U,\G^{(2)}\bigr).
\end{equation}
A simple computation shows that 
\begin{equation}
  K\rmx^\star K^*+\rmx=\p_\l-\p_{\m^\perp}.
\end{equation}

Thus, the transformation of the $\ell$-matrix of
equation~\eqref{eq:cqzejkbvfhqzeuiobvf} by the element $\varsigma\colon
p\mapsto\e^{sp}$ of $\Map_0(U,D)^\l$ (where $D$ is the connected,
simply-connected Lie group with Lie algebra $\dlie$), satisfies:
\begin{equation}
  \left(
    Kj^\star\lcanG{\G^\star,\l,\l^\perp}{j^\star}^*K^*+K\rmx^\star
    K^*+\rmx
  \right)^\varsigma_0=0
\end{equation}
as well as
\begin{equation}
  \left(
    Kj^\star\lcanG{\G^\star,\l,\l^\perp}{j^\star}^*K^*+K\rmx^\star
    K^*+\rmx
  \right)^\varsigma_psp=0.
\end{equation}
Thus, from the uniqueness result of Proposition~\ref{pr:uniqueness}, we
must have:
\begin{equation}
  Kj^\star\lcanG{\G^\star,\l,\l^\perp}{j^\star}^*K^*
  =\left(j\lcanG{\G,\l,\m}j^*\right)^{\varsigma^{-1}}-\p_\l+\p_{\m^\perp}.
\end{equation}



\begin{thebibliography}{10}
\bibitem{AKS00}{A.~Alekseev, Y.~Kosmann-Schwarzbach, Manin pairs and
    moment maps, \emph{J.~Diff.~Geom.}~\textbf{56} (2000), 133--165.} 
\bibitem{AM00}{A.~Alekseev, E.~Meinrenken, The non-commutative Weil
    algebra, \emph{Invent.~Math.}~\textbf{139} (2000), 135--172.} 
\bibitem{D90}{V.~G.~Drinfel{\cprime}d, Quasi-Hopf algebras,
    \emph{Leningrad Math.~J.}~\textbf{1} (1990), 1419--1457.} 
\bibitem{EE03}{B.~Enriquez, P.~Etingof, Quantization of
    Alekseev--Meinrenken dynamical $r$-matrices, in {\it Lie groups and
      symmetric spaces}, 81--98, Amer. Math. Soc. Transl. Ser. 2, 210,
    Amer. Math. Soc. (2003).}
\bibitem{ES01}{P.~Etingof, O.~Schiffmann, On the moduli space of classical
    dynamical $r$-matrices, \emph{Math.~Res.~Lett.}~\textbf{8} (2001),
    157--170.} 
\bibitem{EV98}{P.~Etingof, A.~Varchenko, Geometry and classification of
    solutions of the classical dynamical Yang--Baxter equation,
    \emph{Comm.~Math.~Phys.}~\textbf{192} (1998), 77--120.} 
\bibitem{F95}{G.~Felder, Conformal field theory and integrable systems
    associated to elliptic curves, \emph{Proceedings of the
    International Congress of Mathematicians, Vol.\ 1, 2 (Z\"urich,
    1994)}, Birkh\"auser, Basel (1995), 1247--1255.}
\bibitem{LP05}{L.~C.~Li, S.~Parmentier, On dynamical Poisson groupoids I,
    \emph{Mem.~Amer.~Math.~Soc.}~\textbf{174} (2005), no.~824.}
\bibitem{M87}{K.~C.~H.~Mackenzie, \emph{Lie groupoids and Lie algebroids
      in differential geometry}, London Mathematical Society Lecture
    Note Series, \textbf{124}. Cambridge University Press, Cambridge,
    1987.} 
\bibitem{M99}{K.~C.~H.~Mackenzie, On symplectic double groupoids and the
    duality of Poisson groupoids, \emph{Internat.~J.~Math.}~{\bf 10}
    (1999), no.~4, 435--456.} 
\bibitem{PP05}{S.~Parmentier, R.~Pujol, Quasi-bialgebras and dynamical
    $r$-matrices, \emph{Adv.~Math}~\textbf{197} (2005), 41--85.}
\bibitem{Pthesis}{R.~Pujol, \'Equations de Yang--Baxter dynamiques
    classiques, groupo\"\i{}des de Poisson, quasi-big\`ebres de Lie et
    dualit\'e, Ph.D.\ thesis (in french), Universit\'e Claude Bernard
    Lyon 1 (2005).}
\bibitem{P?}{R.~Pujol, $\A$-dynamical Poisson groupoids, in final
    preparation.}
\bibitem{W88}{A.~Weinstein, Coisotropic calculus and Poisson groupoids,
    \emph{J.~Math.~Soc.~Japan} \textbf{40} (1988), 705--727.}
\end{thebibliography}
\end{document}